\documentclass[1p]{elsarticle}

\usepackage{amsfonts}
\usepackage{amsmath}
\usepackage{amssymb}

\usepackage{abstract}
\usepackage{natbib}
\usepackage{geometry}
\usepackage{multicol}
\usepackage{mathabx}

\usepackage{graphicx}
\usepackage{subfigure,psfrag}

\newcommand{\vs}[1]{\vspace{#1}}

\newcommand{\begitem}{\begin{itemize}}    
\newcommand{\finit}{\end{itemize}}    
\newcommand{\begenum}{\begin{enumerate}}    
\newcommand{\finenum}{\end{enumerate}}    
\newcommand{\begar}{\begin{array}}    
\newcommand{\finar}{\end{array}}    
\newcommand{\begeq}[1]{\begin{equation} \label{#1}}    
\newcommand{\fineq}{\end{equation}}    
\newcommand{\begct}{\begin{center}}    
\newcommand{\finct}{\end{center}}    
\newcommand{\ra}{\rightarrow}

\newcommand{\zun}{\vs{0.1cm}} 
     
\newcommand{\zdeux}{\vs{0.2cm}} 
\newcommand{\ztrois}{\vs{0.3cm}}    
\newcommand{\zcinq}{\vs{0.5cm}} 

\def\tinf{{\rightarrow\infty}}     
\def\tqdninf{{\stackrel{n\tinf}{\longrightarrow}}} 

\def\bE{{\mathbb{E}}}
\def\bN{{\mathbb{N}}}

\def\bP{{\mathbb{P}}}
\def\bI{{\mathbb{I}}}
\def\bV{{\mathbb{V}}}

\def\tinf{{ \rightarrow \infty }}
\def\tqdninf{{ \stackrel{n\tinf}{\longrightarrow} }}

\def\gF{{\gamma_1}}
\def\gG{{\gamma_2}}
\def\gFet{{\gamma_1^*}}
\def\gGet{{\gamma_2^*}}
\def\rF{{\rho_1}}

\def\Fbartn{{\overline F(t_n)}}
\def\unsurFbartn{{\frac 1{\Fbartn} }}
\def\unsurFbartndeux{{\frac 1{\Fbartn^2} }}
\def\Xiet{{X_i^*}}
\def\indicXiet{{\bI_{\Xiet>t_n}}}
\def\logXiet{{\log\left(\frac{\Xiet}{t_n}\right)}}
\def\logxtn{{\log\left(\frac x{t_n}\right)}}
\def\logxtndeux{{\log^2\left(\frac x{t_n}\right)}}
\def\logxtntext{{\log(x/t_n)}}

\def\Fn{{F_n}}
\def\Fbarn{{\overline F_n}}
\def\nb{{\bar{n}}}

\def\moyennesurn{{\frac 1 n \sum_{i=1}^n}}

\def\gamchap{{\widehat\gamma_n}}
\def\gamtilde{{\widetilde\gamma_n}}
\def\xp{{x_{p_n}}}
\def\chapx{{ \hat{x}_{p_n,t_n} }}

\def\Fb{{\overline{F}}}
\def\Gb{{\overline{G}}}
\def\Hb{{\overline{H}}}

\def\Vin{{V_{i,n}}}
\def\Vpin{{V'_{i,n}}}

\def\inttninf{{\int_{t_n}^{\infty}}}
\def\intuninf{{\int_{1}^{\infty}}}
\def\op1{{o_{\bP}(1)}}
\def\opvn12{{o_{\bP}(v_n^{-1/2})}}

\begin{document}
\setlength{\abovedisplayskip}{0.25cm}
\setlength{\belowdisplayskip}{0.25cm}

\newtheorem{theo}{Theorem}
\newtheorem{prop}{Proposition}
\newtheorem{defi}{Definition}
\newtheorem{lem}{Lemma}
\newtheorem{cor}{Corollary}
\newdefinition{rmk}{Remark}
\newdefinition{tb}{Table}

\newcounter{Rem}   
\setcounter{Rem}{-1}    
\refstepcounter{Rem}    
\newcommand{\itemrem}{\addtocounter{Rem}{1} \underline{\rm Remarque \theRem} \ }   
\newcommand{\rem}{\addtocounter{Rem}{1}\hop\underline{\rm Remarque \theRem} \ }   

\begin{frontmatter}

\title{A Lynden-Bell integral estimator for extremes of randomly truncated data} 
 \author{J.~Worms\corref{cor1}\fnref{fn1}}
 \ead{julien.worms@uvsq.fr}
 \author{R.~Worms\fnref{fn2}}
 \ead{rym.worms@u-pec.fr}
 \cortext[cor1]{Corresponding author}
 \fntext[fn1]{Universit\'e de Versailles-Saint-Quentin-En-Yvelines, Laboratoire de Math\'ematiques de Versailles (CNRS UMR 8100), B\^at. Fermat, 45 av. des Etats-Unis, 78035 Versailles, France}
 \fntext[fn2]{UPEMLV, UPEC, Universit\'e Paris-Est, Laboratoire d'Analyse et de Math\'ematiques Appliqu\'ees (CNRS UMR 8050),  F-94010, Cr\'eteil, France}

\date{\today}



\end{frontmatter}

\begin{abstract}
This work deals with the estimation of the extreme value index and extreme quantiles for heavy tailed data, randomly right truncated by another heavy tailed variable. Under mild assumptions and the condition that the truncated variable is less heavy-tailed than the truncating variable, asymptotic normality is proved for both estimators. The proposed estimator of the extreme value index is an adaptation of the Hill estimator, in the natural form of a Lynden-Bell integral. Simulations illustrate the quality of the estimators under a variety of situations. 
\end{abstract}

\section{Introduction}

Extreme value statistics is an active domain of research, with numerous fields of application, and which benefits from an important litterature in the context of i.i.d. data, dependent data, and (more recently) multivariate or spatial data. For univariate data, semiparametric estimation of the tail of the underlying distribution (for instance, estimation of extreme quantiles) requires in the first place accurate estimation of the so-called extreme-value index (e.v.i.).
%
In the recent years, several authors dedicated their efforts to obtaining good estimations of the e.v.i. for incompletely observed data, {\it i.e.} randomly censored or truncated data (note here that, since the interest generally lies in the evaluation of the upper tail of the data, left censoring or left truncation is not a relevant framework, and therefore censoring or truncating are considered from the right). In those contexts, the usual estimators of the e.v.i. need some modifications because otherwise they would lead to erroneous estimations when blindly applied to censored or truncated data. Some references for extreme value estimation in the context of randomly censored observations are \cite{BeirlantDierckxFVGuillou07}, \cite{EinmahlFVGuillou08}, \cite{Worms14}. 
\vspace{0.1cm}

The first published work on extreme values estimation under random truncation was written by L.Gardes and G.Stupfler \cite{GardesStupfler15}, who dealt with heavy-tailed right truncated data (in their work, they provided motivations and many references on main existing results about truncated samples, we refer to \cite{GardesStupfler15} in this regard). The framework of randomly right truncated data will be precisely defined in the next section, let us just sketch it for the moment : we consider $\nb$ independent i.i.d. couples $((X_i,Y_i))_{1\leq i\leq \nb}$ and, among those couples, we only observe those couples which satisfy the condition $X_i\leq Y_i$. The actually observed data will then be noted $((\Xiet,Y_i^*))_{1\leq i\leq n}$. Below, $F$ and $G$ will stand for the respective distributions of $X$ and $Y$, whereas $F^*$ and $G^*$ will stand for the conditional distributions of $X$ and $Y$ given that $X\leq Y$ : the latter two are therefore the distributions of the observed samples $(\Xiet)_{1\leq i\leq n}$ and $(Y^*_i)_{1\leq i\leq n}$. The first objective is to estimate the e.v.i. of $X$.
\vspace{0.1cm}

The original idea in \cite{GardesStupfler15} was to notice that the extreme value indices $\gFet$ and $\gGet$ of $F^*$ and $G^*$ are related by a very simple relation to those of $F$ and $G$, $\gF$ and $\gG$ : they proved that we have indeed (when both $F$ and $G$ are heavy-tailed)
\[
  \gFet = \gF\gG/(\gF+\gG)  \makebox[1.5cm][c]{and}\gGet = \gG.
\] 
These relations readily yield a proposition of estimator for the parameter of interest $\gF$ by relying on usual Hill estimators of $\gFet$ and $\gGet$ :
\begeq{estimGS}
 \hat\gamma_{GS} = \frac{ \hat{\gamma}^*_1(k_1) \hat\gamma_2(k_2)  }{ \hat\gamma_2(k_2) - \hat{\gamma}^*_1(k_1) }    
 \makebox[2.cm][c]{where} 
 \hat{\gamma}^*_1(k_1) = \frac 1{k_1} \sum_{i=1}^{k_1} \log \frac{X^*_{n-i+1,n}}{X^*_{n-k_1,n}}
 \makebox[1.5cm][c]{and}
 \hat\gamma_2(k_2) = \frac 1{k_2} \sum_{i=1}^{k_2} \log \frac{Y^*_{n-i+1,n}}{Y^*_{n-k_2,n}}
\fineq 
where $X^*_{1,n}\leq \ldots\leq X^*_{n,n}$ and $Y^*_{1,n}\leq \ldots\leq Y^*_{n,n}$ denote the usual order statistics of both samples, and $k_1$ and $k_2$ are the number of upper observations which are kept for estimating $\gFet$ and $\gGet$. 
\vspace{0.1cm}

The authors of \cite{GardesStupfler15} also investigated the behavior of an estimator of $F$ in the upper tail, and therefore provided a Weissman-type estimator of extreme quantiles in this truncation context and proved its asymptotic normality. However, their results suffer from some kind of calibration problem, since they are proved only under the condition that one of the numbers $k_1$ and $k_2$ of order statistics used for estimating $\gFet$ and $\gG$ must grow to infinity faster than the other. The question of getting rid of this restriction was addressed in the prepublication \cite{NecirTroncatureHill15}.
\vspace{0.1cm}

In this work, we consider the same framework of randomly right-truncated heavy-tailed data, but adopt a new method for defining an estimator of the extreme value index $\gF$ of the truncated sample : in Section 2, this estimator $\gamchap$ is defined as some Lynden-Bell integral, requiring a single threshold to be chosen, and asymptotic normality is proved for $\gamchap$ as well as for an estimator  of extreme quantiles, under appropriate but mild conditions. Section 3 is devoted to a simulation study illustrating the performance of the defined estimators (with a tentative comparison to the performance of the estimator defined in \cite{GardesStupfler15}), and Sections 4 and 5 respectively contain a conclusion and  the proofs of the results. The appendix recalls important (and needed) results, previously published in the litterature, and contains as well a technical lemma which is repeatedly used in the proofs section.

\section{Framework and statement of the results}

\subsection{Notations and definition of the estimators}

Let $((X_i,Y_i))_{1\leq i \leq \nb}$ be $\nb$ independent copies of a couple $(X,Y)$, where $X$ and $Y$ are positive independent random variables having respective cumulative distribution functions  $F$ and $G$. For convenience, we suppose that the lower endpoints of $F$ and $G$ are both equal to $0$ (but this will have no influence on the results, since only the highest data values are retained for tail estimation). We assume in this work that $X$ and $Y$ are heavy-tailed distributed, meaning that $1-F$ and $1-G$ (also assumed to be continuous) are regularly varying with respective indices $-1/\gF$ and $-1/\gG$ where $\gF$ and $\gG$ are $>0$.  
\zdeux

We only observe the couples $(X_i,Y_i)$ which satisfy $X_i\leq Y_i$ : in other words, the original data $X_i$ are randomly truncated from the right by the $Y_i$, and the actually observed sample is $((X_i^*,Y_i^*))_{1\leq i \leq N}$, where $N$ follows the ${\cal B}(\nb,p)$ distribution, $p$ denoting the (unknown) probability of non-truncation $p=\bP(X\leq Y)$. Consequently, the distribution of the $\Xiet$ becomes
\begeq{relationFetF}
 F^*(x) = \bP( X\leq x | X\leq Y ) = \frac 1 p \int_0^x \Gb(t) dF(t).
\fineq
Conditionally on $N=n$, the couples $(X^*_1,Y^*_1),\ldots,(X^*_N,Y^*_N)$ are independent and identically distributed, and $X^*_i$ is no longer independent of $Y^*_i$. It is important to note that, in the sequel, we will work conditionnaly on $N=n$, where $n$ is some deterministic sample size, and we will therefore handle the sample  $(X^*_1,Y^*_1),\ldots,(X^*_n,Y^*_n)$ without further reference to $N$.
\zdeux

In this work, $F_n$ will denote the classical Lynden-Bell (nonparametric maximum likelihood) estimator of $F$, namely
\[
 F_n(x) = \prod_{\Xiet>x} \left( 1 - \frac 1{nC_n(\Xiet)} \right)
 \makebox[2.cm][c]{where}
 C_n(x) = \moyennesurn \bI_{\Xiet \leq x\leq Y_i^*}
\]
(with the usual convention that a product on the empty set equals $1$), where $C_n$ is the estimator of the function $C$
\begeq{defC}
C(x)=\bP(X\leq x\leq Y | X\leq Y) = p^{-1}\bar G(x) F(x)
\fineq
 which plays an important role in the analysis of truncated data. Note that $F_n$ is very close to, but different strictly speaking, from the estimator of $F$ considered in \cite{GardesStupfler15} ($F_n$ takes rational values, which is not the case of the latter).
\ztrois\\

Our goal  is to adapt  the famous Hill estimator in the context of right-truncation. It is well known that (see Remark 1.2.3 in \cite{deHaanFerreira2006} for instance)
\[
 \bE \, [\log(X/t) \, | \, X>t  ] =  \frac{1}{\Fb(t)} \int_t^{\infty} \log (x/t) \ dF(x)
\]
tends to $\gF$ as $t\ra +\infty$.  If $(t_n)$  is  a sequence of positive thresholds growing to infinity with $n$,  we can then define a random version of $\phi(x) = (\Fb(t))^{-1}   \log(x/t)\bI_{x>t}$ by $\hat{\phi}_n(x) = (\Fb(t_n))^{-1}  \log(x/t_n)\bI_{x>t_n}$ 
 and consequently,  a natural adaptation of the Hill estimator for $\gF$ is  (see relations (1.9) and (1.10) in \cite{StuteWang2008}, in the left-truncation case, for details about Lynden-Bell integrals)
\[
\gamchap = \int \hat{\phi}_n(x) dF_n(x) = \frac 1 n \sum_{i=1}^n   \hat{\phi}_n(\Xiet)\frac{\Fn(\Xiet)}{C_n(\Xiet)}, 
\] 
which leads to 
\begeq{defgamchap}
 \gamchap = \frac 1  {n\Fbarn(t_n)} \sum_{i=1}^{n} \logXiet \frac{\Fn(\Xiet)}{C_n(\Xiet)} \indicXiet
\fineq
\zun 

Note that this principle has already been successfully applied in the censoring framework in \cite{Worms14} (see equation $(7)$), where the role of  Lynden-Bell estimator was played by the Kaplan-Meier estimator. However, here, the threshold $t_n$ is deterministic instead of being an order statistic. 
The asymptotic properties of $\gamchap$ are stated in Theorem \ref{TLCgamchap}. Naturally, the lighter the truncation, the closer our estimator $\gamchap$ gets to the usual Hill estimator. (?)
\zdeux

We will use this estimator of the tail index $\gF$ in order to estimate an extreme quantile, following a classical scheme. More precisely, let $(p_n)$ be some sequence of quantiles orders tending to $0$, such that $p_n=o(\Fb(t_n))$. If $x_{p_n}$ denote the quantile of $F$ of order $1-p_n$, {\it i.e.} solving $\Fb(x_{p_n})=p_n$, then, in this heavy tailed context (see (\ref{Ordre1}) below), it is easy to see that we can define an estimator $\chapx$ of $\xp$ as 
\begeq{defchapx}
\chapx = t_n \left( \frac{\Fbarn(t_n)}{p_n} \right)^{\gamchap}.
\fineq
In the situation of untruncated data, this is a classical estimator for an extreme quantile based on the approximation of the log relative excesses by a Pareto distribution in the heavy-tailed context, where $F_n$ is in this case the empirical distribution function.

\subsection{Assumptions and results}

The first order condition assumed in this work is the following 
\begeq{Ordre1}
\Fb \in RV_{-1 / \gF} \mbox{ \ and \ } \Gb \in RV_{-1 / \gG} \mbox{ \ with \ } 0 < \gF < \gG. 
\fineq
In other words, we assume that the tail of the truncating variable $Y$ is heavier than the tail of the variable $X$ of interest. This condition is needed in many occasions in the proofs of our results, and is due to the presence (in (\ref{defgamchap})) of the Lynden Bell estimator, evaluated in the tail. Note that this implies the finiteness of the integral $\int_0^{\infty} dF(x)/\Gb(x)$ (which is a sufficient condition sometimes stated in papers dealing with the asymptotic normality of $F_n$). 
\zun

Moreover, if we note  $l_F$ the slowly varying function associated to $F$ ({\it i.e.} such that  $\Fb (x)= x^{-1/\gF} l_F(x)$), the second order condition we consider is the classical $SR2$ condition for $l_F$ (see \cite{BinghamGoldieTeugels}),
\begeq{Ordre2}
\forall x >0 , \ \frac{l_F(tx)}{l_F(t)} -1 \ \stackrel{t\tinf}{\sim} \ h_{\rF}(x) \  g(t) \hspace{0.3cm} (\forall x >1)
\fineq
where $g$ is a positive mesurable function, slowly varying with index $\rF$, and $h_{\rF}(x) = \frac{x^{\rF}-1}{\rF}$  when  $\rF <0$, or  $h_{\rF}(x) =\log x$ when  $\rF = 0$.
\zun

The first assumption on the threshold sequence $(t_n)$ will be that,  if we note  $\Hb = \Fb \ \Gb$ (note that $H$ is the distribution function of $\min(X,Y)$), $(t_n)$ satisfies
\begeq{nHbar}
n \Hb (t_n) \ \tqdninf \ + \infty. 
\fineq
The asymptotic normality result will then require the following condition on $(t_n)$ : 
\begeq{vng}
\sqrt{n \Hb (t_n)} g(t_n) \ \tqdninf \ \lambda \; \mbox{ \ for some } \lambda >0. 
\fineq

\begin{theo}\label{TLCgamchap}
Under assumptions $(\ref{Ordre1})$, $(\ref{Ordre2})$,  $(\ref{nHbar})$ and  $(\ref{vng})$, as $n$ tends to infinity, 
\[
\sqrt{n \Hb (t_n)}  (\gamchap - \gF)     \stackrel{\cal L}{\longrightarrow}  {\cal N} \left(\lambda m , s^2 \right),
\]
where $m = \left\{ \begar{ll}    \frac{\gF^2}{1-\gF \rF} & \mbox{ if } \rF <0, \zdeux\\  \gF^2 & \mbox{ if } \rF = 0. \finar \right. $ 
\ and  \ 
$\displaystyle s^2=p \gF^2 \left(1+\left(\frac{\gF}{\gG} \right)^2\right) \left(1-\frac{\gF}{\gG} \right)^{-3}$. 
\end{theo}

Let us now turn to the results about the extreme quantile estimator defined in (\ref{defchapx}). Suppose that the sequence of quantile orders $(p_n)$, tending to $0$, satisfies the condition 
\begeq{condsurpn}
 \Fbartn / p_n \ \tqdninf \ +\infty.
\fineq

\begin{theo}\label{TLCchapx}
Under (\ref{condsurpn}) and the assumptions of Theorem \ref{TLCgamchap}, setting \ $d_n=\Fbartn/p_n$,  \ if $\rF <0$ and
\begeq{condn}
\sqrt{n \Hb(t_n)}\, \big/ \log d_n \rightarrow \infty,
\fineq
as $n$ tends to $\infty$ then
\[
\frac{\sqrt{n \Hb(t_n)}}{\log d_n} \left( \frac{\chapx}{\xp} - 1\right) \stackrel{\cal L}{\longrightarrow}  {\cal N} \left(\lambda m , s^2 \right)
\]  
\end{theo}

\section{Finite sample behaviour}

In this section, we illustrate our results by presenting some graphics (issued from an extensive study) corresponding to the comparison, in terms of bias and root mean squared error (RMSE), of our new estimator $\gamchap$ (defined in (\ref{defgamchap})) with the existing estimator $\hat \gamma _{GS}$ (defined in equation (\ref{estimGS})) issued from \cite{GardesStupfler15}, for two classes of heavy-tailed distributions:
\begitem
\item Burr$(\beta,\tau,\lambda)$ with distribution function $1-(\frac{\beta}{\beta+x^{\tau}})^{\lambda}$, for which the e.v.i. is $\frac{1}{\lambda \tau}$.
\item Frechet$(\gamma)$ with distribution function $\exp(-x^{-1/\gamma})$, for which the e.v.i. is $\gamma$.
\finit

\begin{figure}[htp]
\centering
\subfigure[Burr$(10,4,1)$ truncated by Burr$(10,2,1)$]{
\label{BurrBurrforte}
\includegraphics[height=3.0cm,width=.46\textwidth]{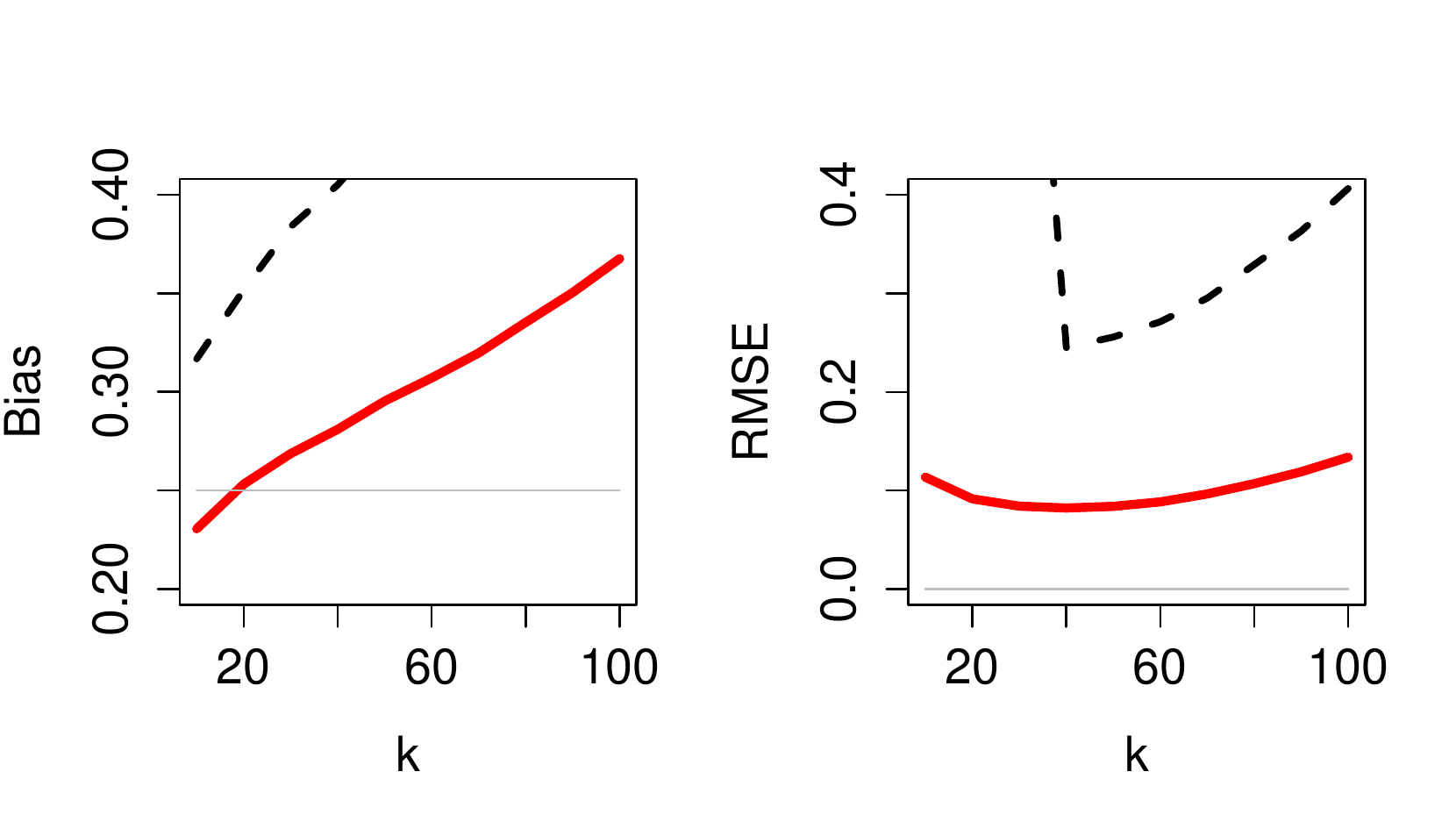}}
\hspace{.1in}
\subfigure[Burr$(10,4,1)$ truncated by Burr$(10,1,1/2)$]{
\label{BurrBurrfaible}
\includegraphics[height=3.0cm,width=.46\textwidth]{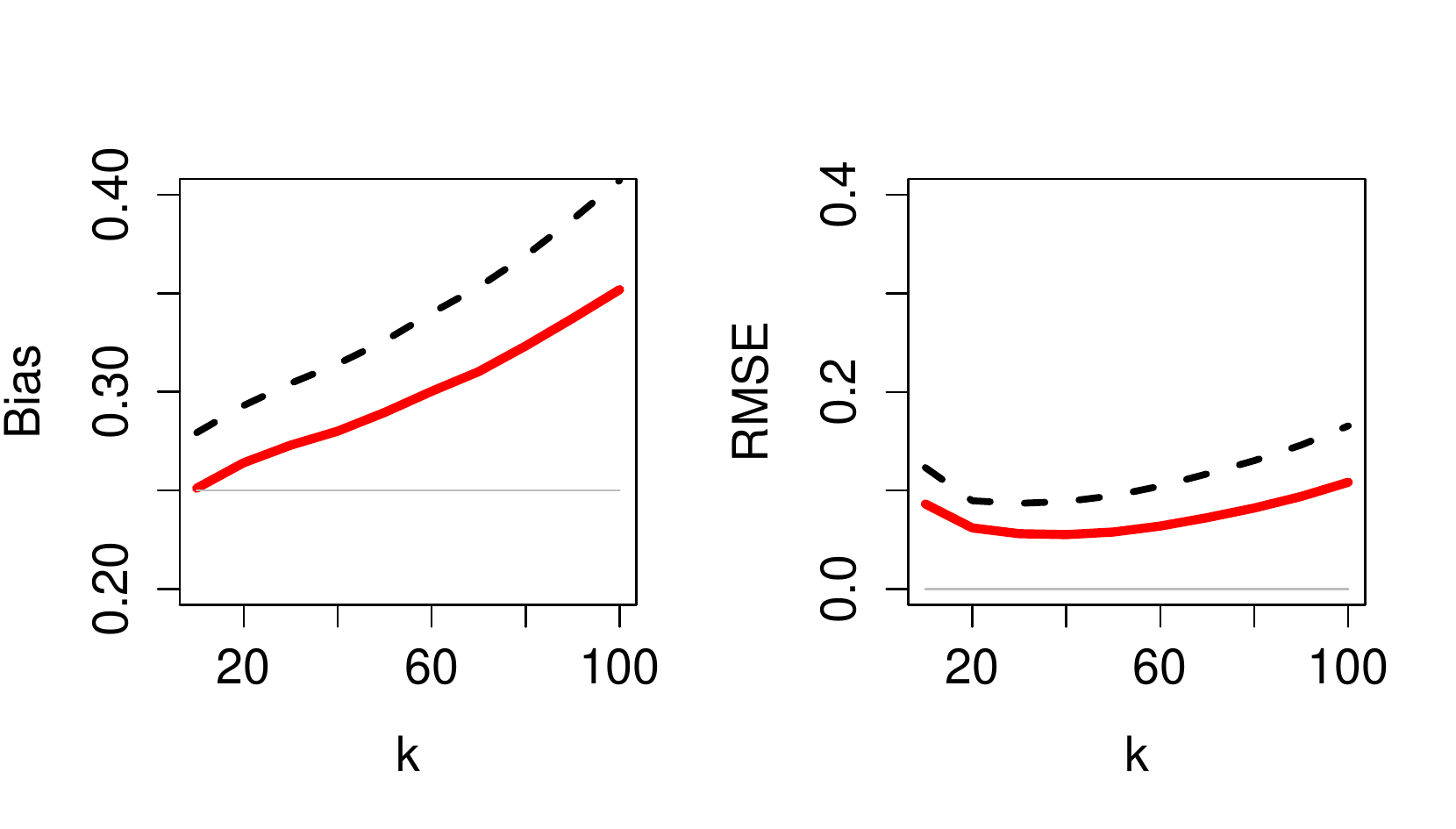}}
\subfigure[Frechet$(1/4)$ truncated by Frechet$(1/2)$]{
\label{FrechetFrechetforte}
\includegraphics[height=3.0cm,width=.46\textwidth]{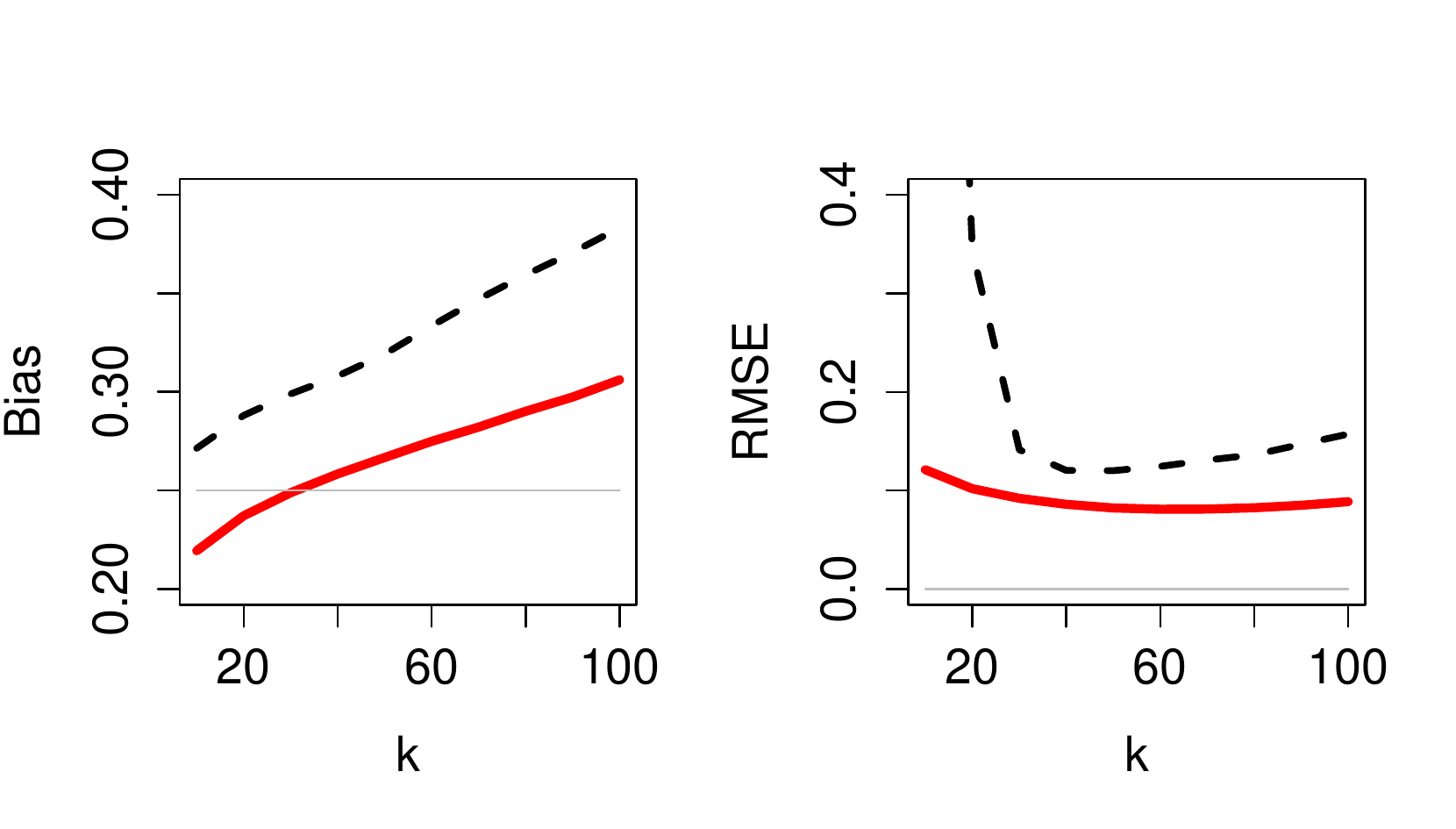}}
\hspace{.1in}
\subfigure[Frechet$(1/4)$ truncated by Frechet$(2)$]{
\label{FrechetFrechetfaible}
\includegraphics[height=3.0cm,width=.46\textwidth]{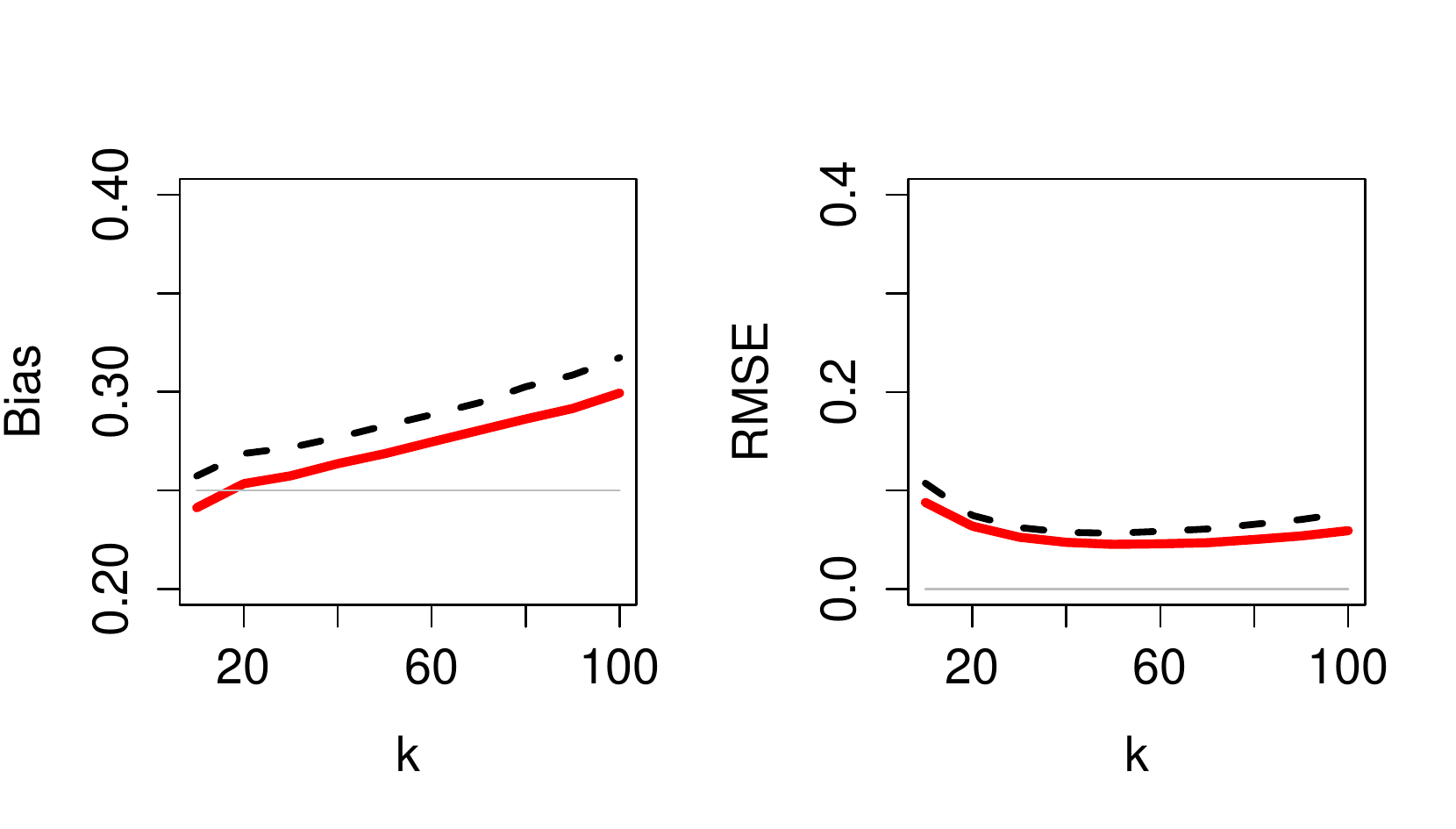}}
\subfigure[Burr$(10,4,1)$ truncated by Frechet$(1/2)$]{
\label{FrechetFrechetforte}
\includegraphics[height=3.0cm,width=.46\textwidth]{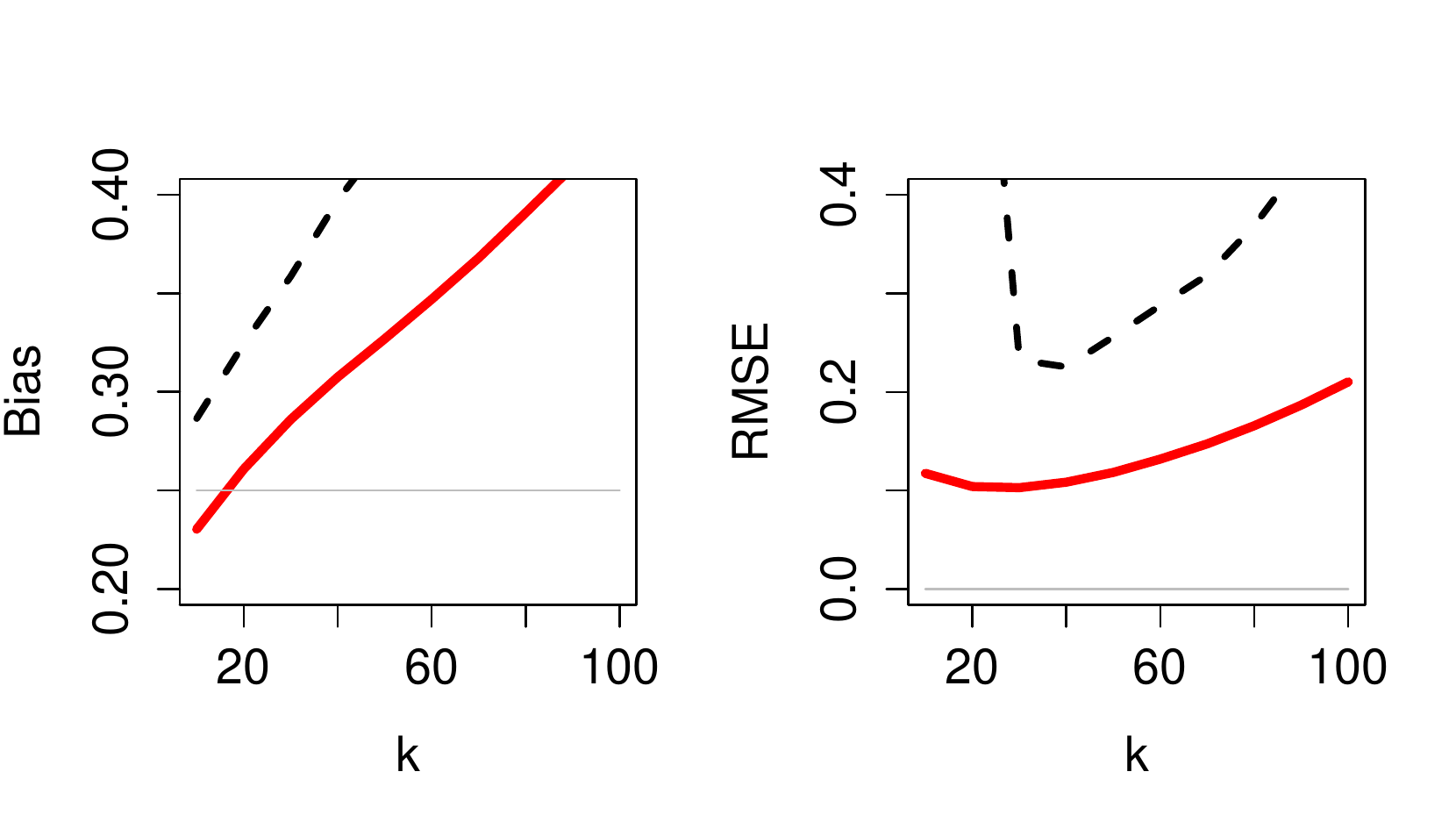}}
\hspace{.1in}
\subfigure[Burr$(10,4,1)$ truncated by Frechet$(2)$]{
\label{FrechetFrechetfaible}
\includegraphics[height=3.0cm,width=.46\textwidth]{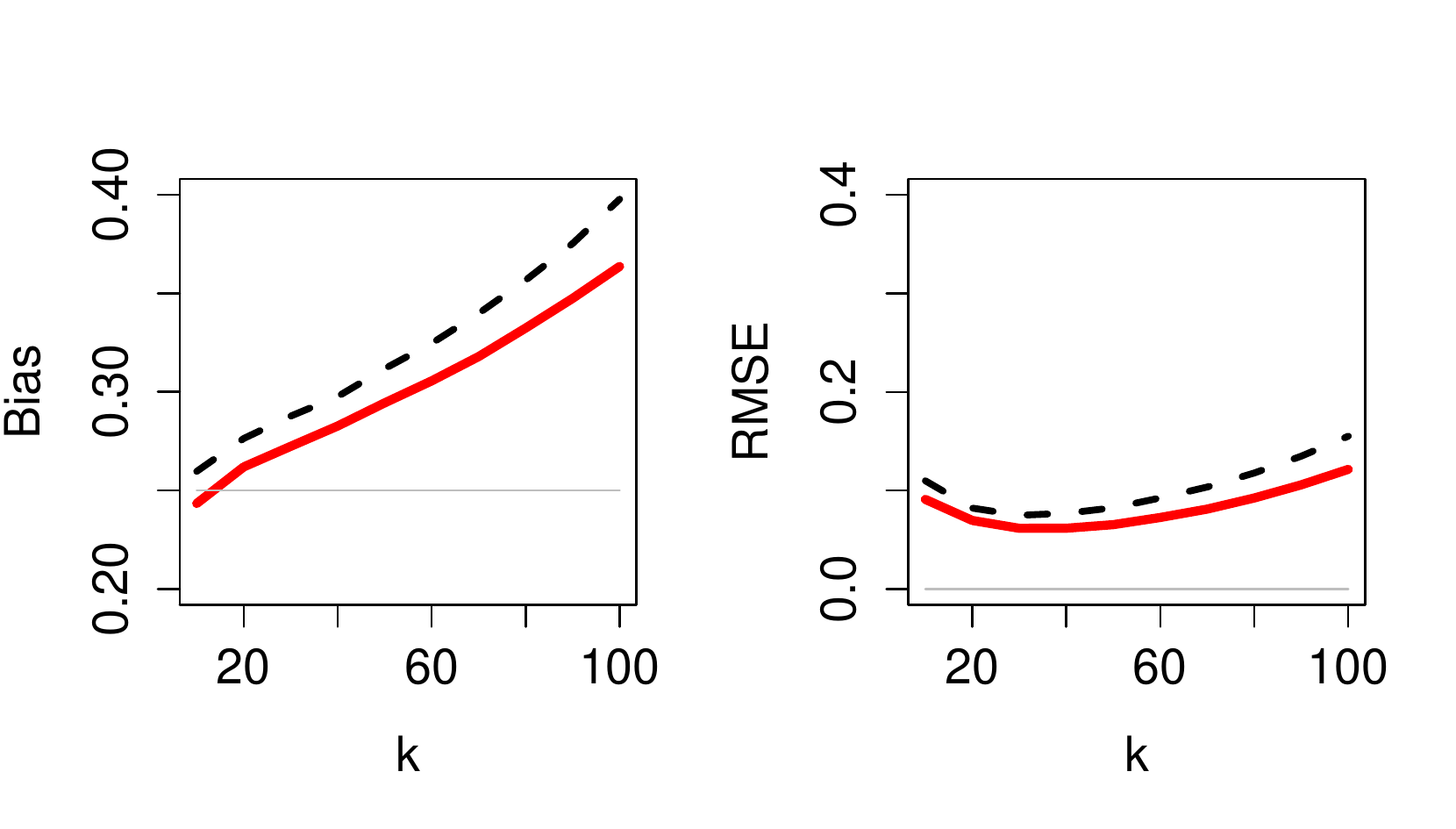}}
\caption{Comparison of bias  and RMSE (respectively left and right in each subfigure) for $\gamchap$ (plain) and $\hat\gamma _{GS}$ (dashed) where $\gF =1/4$, $\gG=1/2$ and $\alpha=2/3$ (important truncation) for subfigures (a),(c),(e), and where $\gF =1/4$, $\gG=2$ and $\alpha=8/9$ (mild truncation) for subfigures (b),(d),(f)}
\label{figgamchap}
\end{figure}

Note that, in those simulations, we used the random threshold $X_{n-k_n,n}^*$ (where $1\leq k_n<n$) instead of a deterministic threshold $t_n$ in the definition of $\gamchap$, and we also considered $k_1=k_2$ in the definition of $\hat \gamma _{GS}$, which is out of the scope of Theorem 3 in \cite{GardesStupfler15} (but the authors themselves restricted their simulations to this situation, which was then presented as more manageable and convenient). Note that making $n$ vary did not provide notable findings, so we kept the number $n$ of actual observation fixed.
\zdeux

We simulated $2000$ random samples of size $n=200$ in 6 different situations : 3 choices of families of distributions (Burr truncated by another Burr, Fr\'echet truncated by another Fr\'echet,  and Burr truncated by a Fr\'echet) combined with 2 choices of truncation strength. This strength is measured by the ultimate probability  $\alpha:= \frac{\gG}{\gF+\gG}$ of non-truncation in the tail (for a proof of this formula, see \cite{EinmahlFVGuillou08}), which is distinct from the overall $p=\bP(X\leq Y)$ : two values were considered, $\alpha=2/3$ (for $\gF=1/4$ and $\gG=1/2$, {\it i.e.} important truncation) and $\alpha=8/9$ (for $\gF=1/4$ and $\gG=2$, {\it i.e.} mild truncation). The results are contained in Figure \ref{figgamchap}, where bias and RMSE are plotted against different values of $k_n$, the number of excesses used.
\zdeux

\begin{figure}[htp]
\centering
\subfigure[Burr$(10,4,1)$ truncated by Burr$(10,2,1)$]{
\label{BurrBurrforte}
\includegraphics[height=3.0cm,width=.46\textwidth]{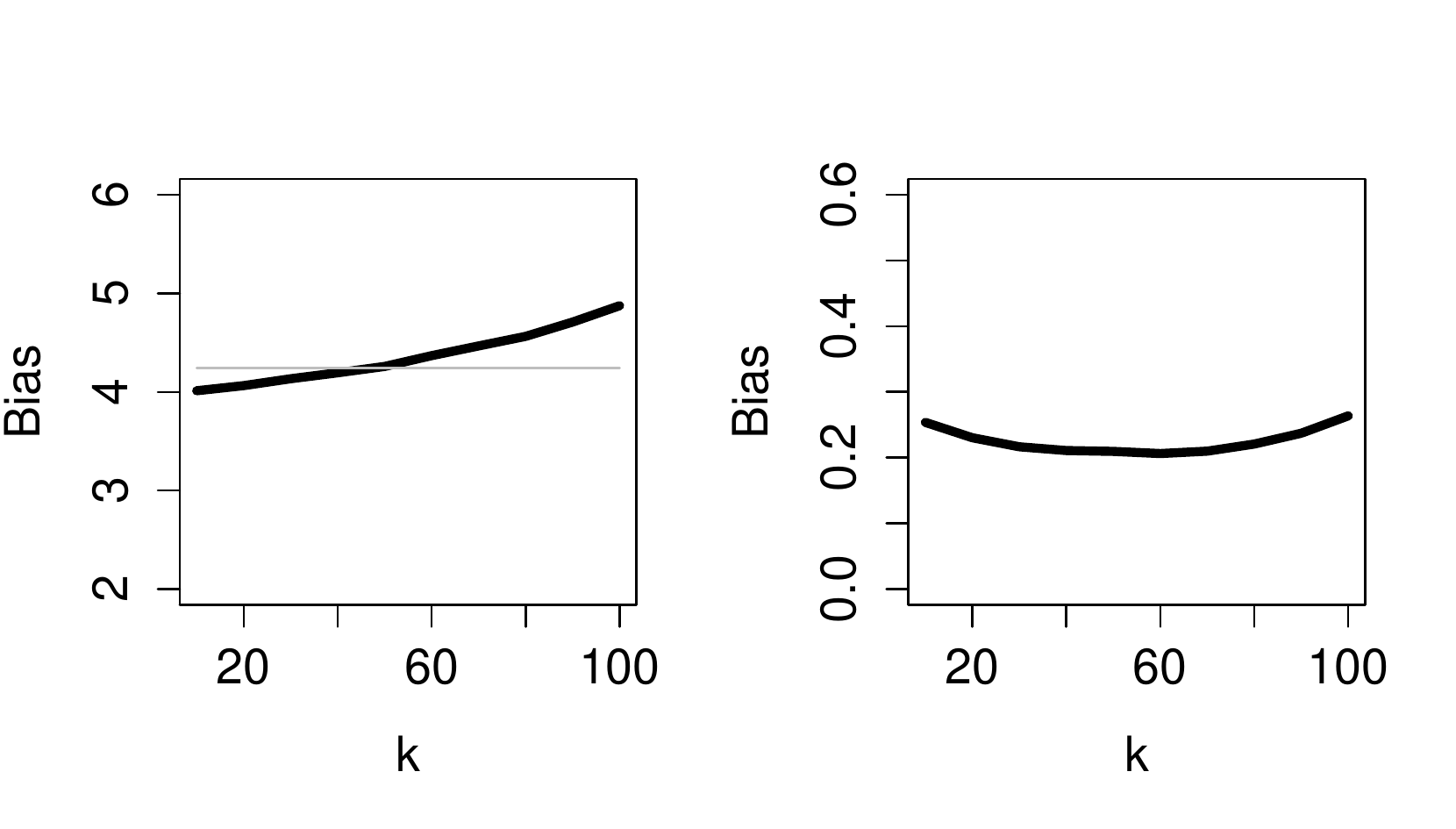}}
\hspace{.1in}
\subfigure[Burr$(10,4,1)$ truncated by Burr$(10,1,1/2)$]{
\label{BurrBurrfaible}
\includegraphics[height=3.0cm,width=.46\textwidth]{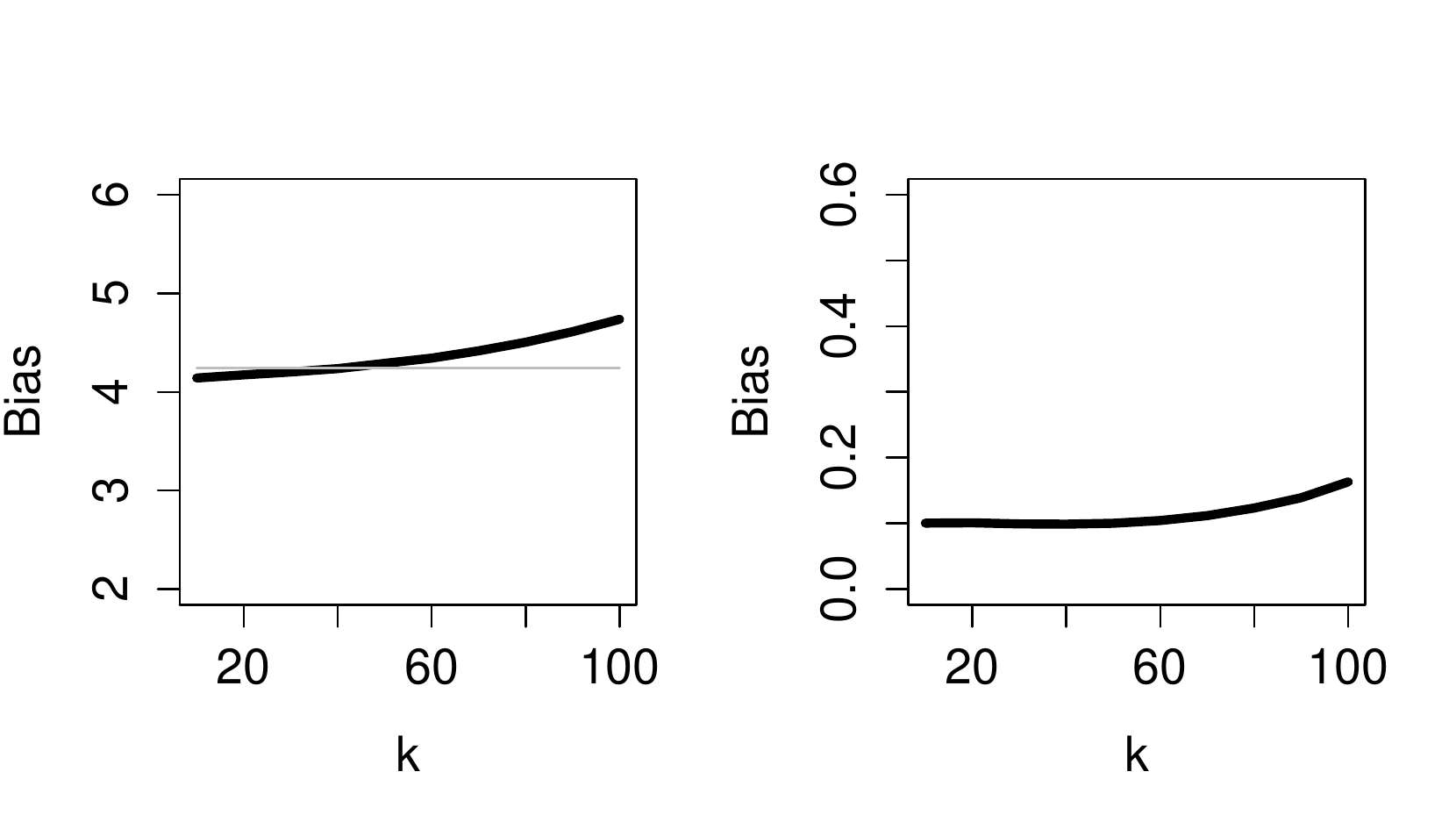}}
\subfigure[Frechet$(1/4)$ truncated by Frechet$(1/2)$]{
\label{FrechetFrechetforte}
\includegraphics[height=3.0cm,width=.46\textwidth]{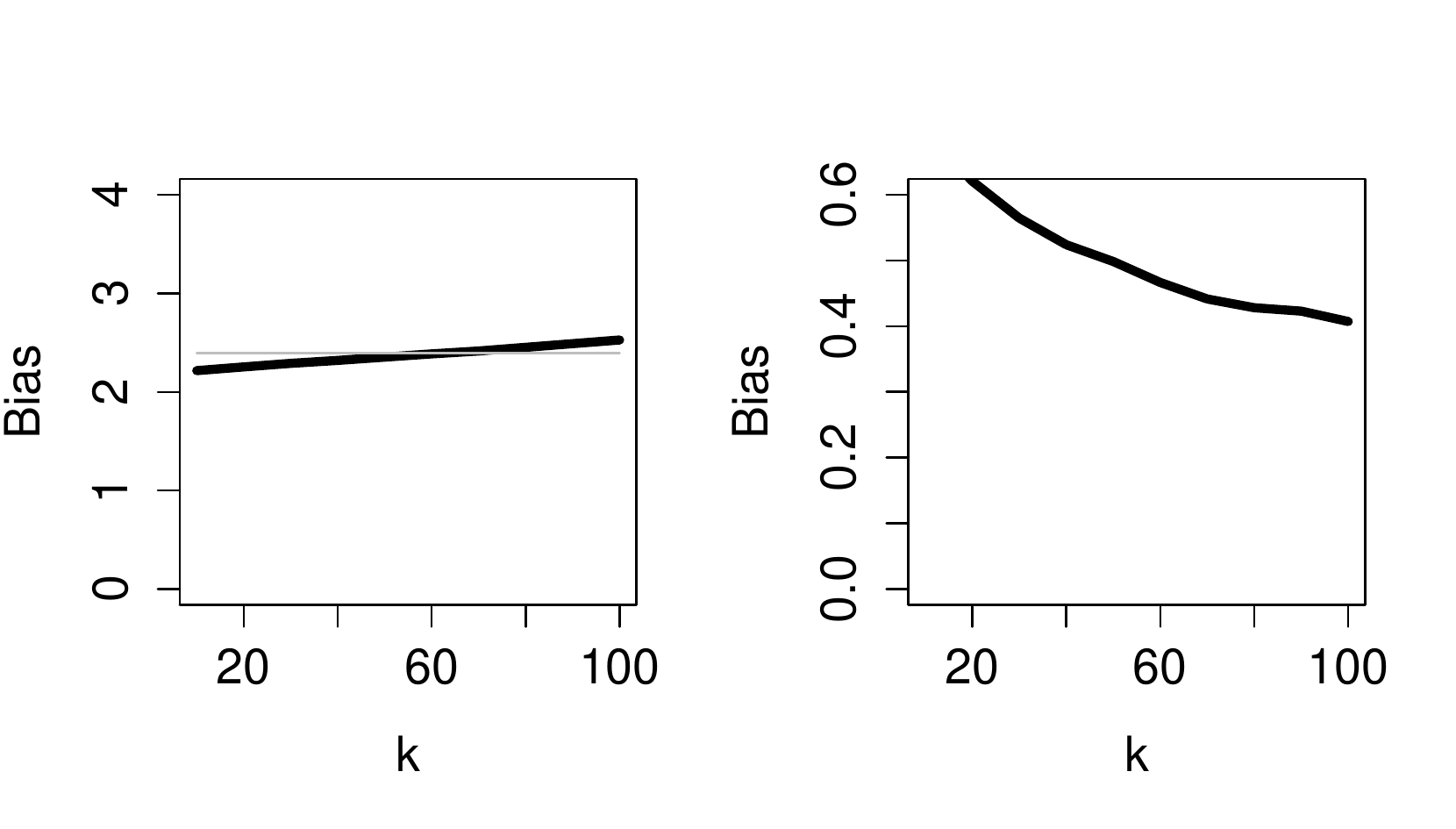}}
\hspace{.1in}
\subfigure[Frechet$(1/4)$ truncated by Frechet$(2)$]{
\label{FrechetFrechetfaible}
\includegraphics[height=3.0cm,width=.46\textwidth]{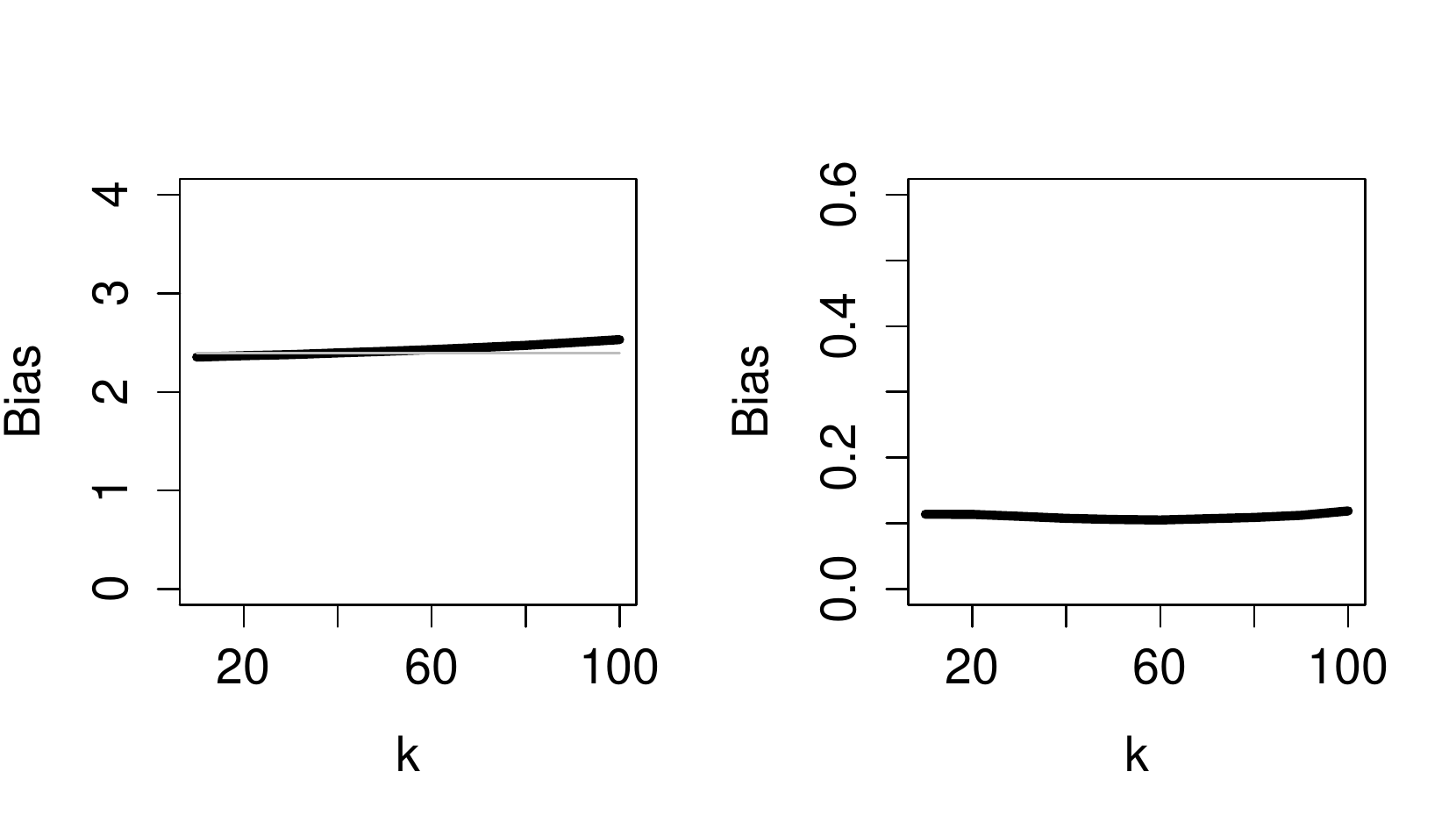}}
\subfigure[Burr$(10,4,1)$ truncated by Frechet$(1/2)$]{
\label{FrechetFrechetforte}
\includegraphics[height=3.0cm,width=.46\textwidth]{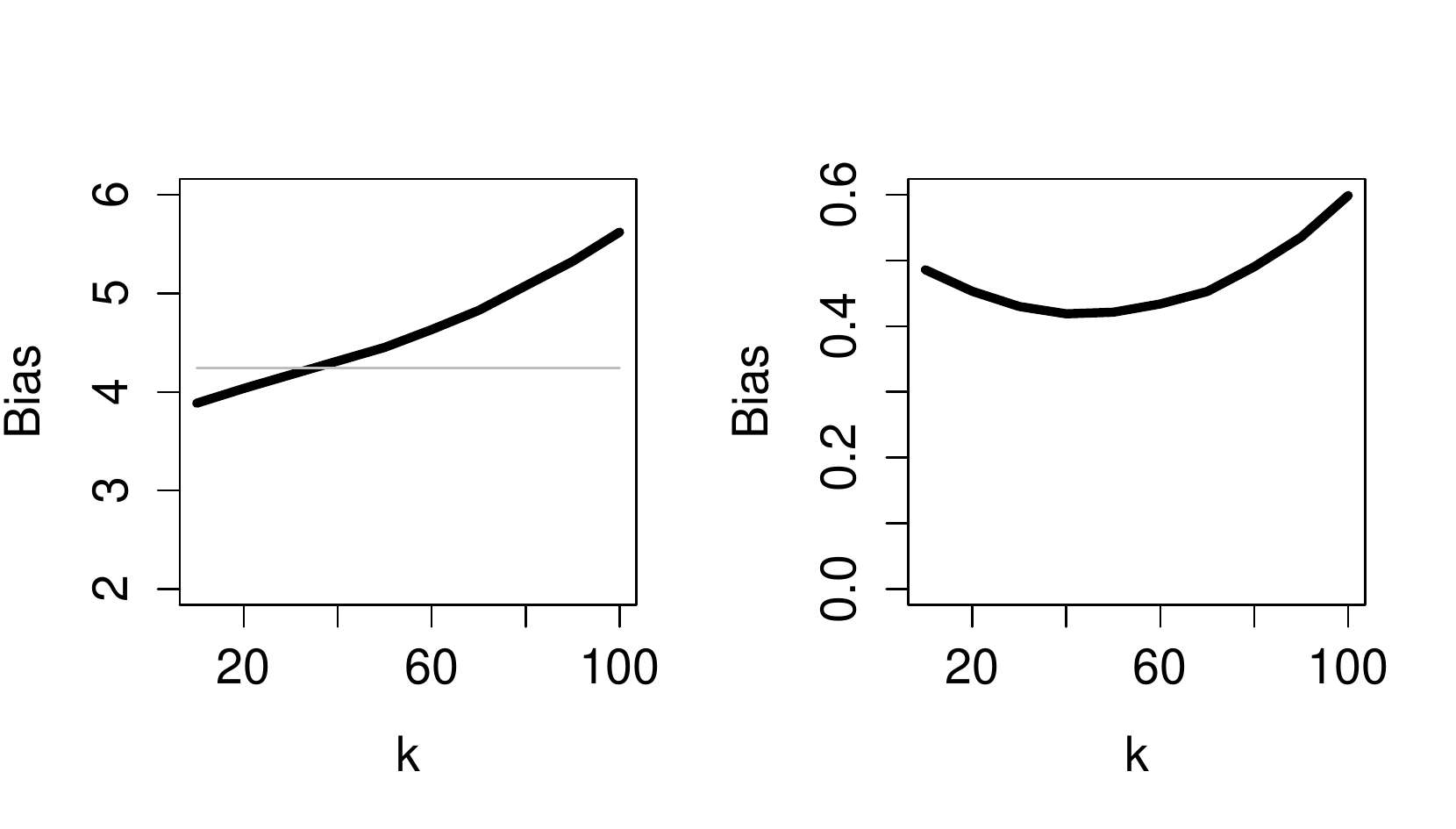}}
\hspace{.1in}
\subfigure[Burr$(10,4,1)$ truncated by Frechet$(2)$]{
\label{FrechetFrechetfaible}
\includegraphics[height=3.0cm,width=.46\textwidth]{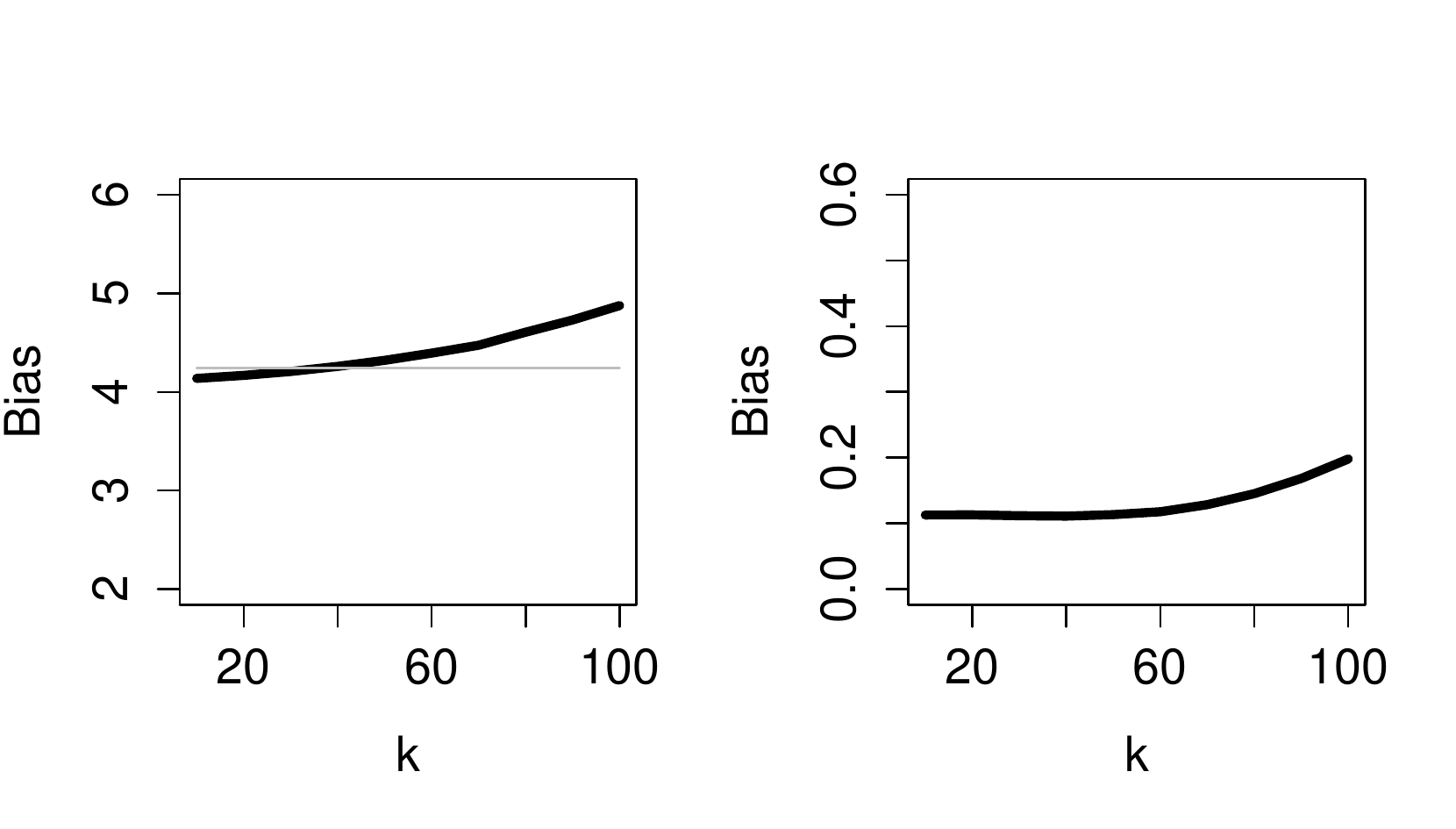}}
\caption{Bias  and RMSE (respectively left and right in each subfigure) for $\chapx$ where $\gF =1/4$, $\gG=1/2$ and $\alpha=2/3$ (important truncation) for subfigures (a),(c),(e), and where $\gF =1/4$, $\gG=2$ and $\alpha=8/9$ (mild truncation) for subfigures (b),(d),(f)}
\label{figxchap}
\end{figure}

This section also contains graphics illustrating the behaviour of our extreme quantile estimator $\chapx$ of $\xp$ (again computed with the random threshold $X_{n-k_n,n}^*$ instead of $(t_n)$. Under the same simulation framework described above, we considered the estimation of the extreme quantile $\xp$ with $p_n=0,03$. Results are displayed in Figure \ref{figxchap}.
\zdeux

The main conclusion we can deduce from our intensive simulation study is that our estimator $\gamchap$ seems to behave systematically better (both in terms of bias and RMSE) than the existing estimator $\hat \gamma _{GS}$ used with $k_1=k_2$, whatever the distributions and the value of $\alpha$ are (and changing the sample size yields the same conclusion). Nonetheless, the comparison may be a bit delicate since the properties of $\hat\gamma_{GS}$ are only proved when the two numbers $k_1$ and $k_2$ are quite distant from each other. On the other hand, the performance of our estimator clearly diminishes when the ultimate proportion of non-truncation $\alpha$ decreases (which is equivalent to $\gF$ getting closer to $\gG$, which notably increases the asymptotic variance of our estimator) but this phenomenon also holds (and to a greater extent) for $\hat \gamma _{GS}$. According to our investigations, and unsurprisingly, a small value of $\rF$ also implies a lesser performance. 
And concerning the bias, since our estimator of $\gF$ is based on the same idea as the Hill estimator in the complete data setting, the relatively high bias observed is neither surprising nor unbearable ; and it is always lower than the bias of $\hat \gamma _{GS}$.

Concerning our new extreme quantile estimator $\chapx$, the finite sample behaviour seems quite satisfying, even if its performances depend on the value of $p_n$ and of the truncation strength. 

\section{Conclusion}

This paper addressed the problem of estimating tails (extreme value index $\gF$ and extreme quantiles) of randomly right-truncated data, when both the truncated and the truncating variables are heavy-tailed. This framework was first considered in \cite{GardesStupfler15}, where a first proposition 
of estimator of $\gF$ was provided. We propose here an alternative approach, leading to an estimator of $\gF$ which takes the form of a Lynden-Bell integral of some particular function, and is therefore a sort of natural version of the Hill estimator in this truncation context. Contrary to the situation of \cite{GardesStupfler15} (for which the choice of the numbers of upper order statistics $k_1$ and $k_2$ in the estimator $\hat\gamma_{GS}$ defined in (\ref{estimGS}) could remain very delicate in practice), a single tuning parameter has to be determined (the threshold $t_n$, or in practice the number of upper order statistics), and experimental results are very encouraging. 
\zun

Concerning the asymptotic normality result for our estimator,  the restriction that the truncating variable has a heavier tail than the truncated variable seems to be unavoidable, and improving the performance in term of bias is an open problem, as is the extension of the approach to truncated data with non-negative extreme value index.

\section{Proofs of the results}

\subsection{Proof of Theorem \ref{TLCgamchap}} 

We introduce the following important notations : first
\begeq{defgamtildetche}
 \gamtilde = \moyennesurn \Vin
 \makebox[2.cm][c]{ where } \Vin = \unsurFbartn \logXiet \frac{F(\Xiet)}{C(\Xiet)} \indicXiet  
\fineq
The variables $\Vin$ are independent and identically distributed and, using (\ref{relationFetF}),
we readily have $\bE(V_{1,n}) = \unsurFbartn \inttninf \logxtntext dF(x)$, which converges to $\gF$. Then we consider two (very close but different anyway) estimators of the cumulative hazard function $\Lambda$  of $X$, $\Lambda=-\log F$ : for any $t$, let (for the first definition below, $F_n(t)$ is supposed $>0$ though)
\begeq{defLambdan}
 \Lambda_n(t) = -\log \Fn(t)  \makebox[2.cm][c]{ and } 
 \hat\Lambda_n(t) = \sum_{\Xiet>t} \frac 1 {nC_n(\Xiet)} . 
\fineq
We will later approach $\hat\Lambda_n(t_n)/\Fbartn$ by $\moyennesurn \Vpin$, where the i.i.d. variables $\Vpin$ are defined by
\begeq{defVpin}
 \Vpin = \frac { \indicXiet } { \Fbartn C(\Xiet) }   \makebox[2.cm][c]{ with }  \bE(V'_{1,n}) = \frac{\Lambda(t_n)}{\Fbartn}.
\zdeux
\fineq
Finally we set $W_{i,n}=\Vin - \bE(V_{1,n})$ and $W'_{i,n}=\Vpin - \bE(V'_{1,n})$, as well as
\[
 \Delta_n = \Fbarn(t_n)/\Fbartn    \makebox[1.5cm][c]{ and }    v_n=n\Hb(t_n)
\]  
\zun\\
Before proceeding to the proof of Theorem \ref{TLCgamchap}, let us state some lemmas (compl\`eter bien s\^ur les conditions/hypoth\`eses...)

\begin{lem}\label{lemmegamchapgamtilde}
 Under condition $(\ref{Ordre1})$, we have 
 $ \Delta_n \gamchap - \gamtilde = o_{\bP}(v_n^{-1/2})  $.
\end{lem} 

\begin{lem}\label{lemmeDeltan}
 Under conditions $(\ref{nHbar})$ and $(\ref{Ordre1})$, the sequence $(\Delta_n)$ converges to $1$ in probability.
\end{lem} 

\begin{lem}\label{lemmeContigLambdas}
 If $T= \max\{\Xiet ; nC_n(\Xiet) =1\}$ and $A_n=\{T\leq t_n\}$, then, under condition $(\ref{Ordre1})$, we have 
 $$ 
 \frac{\sqrt{v_n}}{\Fbartn} \bI_{A_n} (\Lambda_n(t_n) - \hat\Lambda_n(t_n)) = o_{\bP}(1)  .
 $$
\end{lem} 

\begin{lem}\label{lemmeTLClambdan}
 Under conditions $(\ref{nHbar})$ and $(\ref{Ordre1})$,
\begeq{passageLambdanWn}
 \sqrt{v_n} \frac{ \hat\Lambda_n(t_n) - \Lambda(t_n)} { \Fbartn } = 
 \sqrt{v_n} \overline{W}'_n  +  o_{\bP}(1)  .
\fineq
\end{lem} 

For the next two lemmas, note that quantities $s^2$ and $m$ have been defined in the statement of Theorem \ref{TLCgamchap}).

\begin{lem}\label{lemmeWnWtildenTLC}
 Under conditions $(\ref{nHbar})$ and $(\ref{Ordre1})$,  the sequences $\sqrt{v_n} \overline{W}_n$, 
 $\sqrt{v_n} \overline{W}'_n$ and $\sqrt{v_n} (\overline{W}_n - \gF\overline{W}'_n)$ converge in distribution to centered gaussian distributions of respective variances $2p\gF^2/(1-\gF/\gG)^3$, \ $p/(1-\gF/\gG)$ and $s^2$.
\end{lem} 

\begin{lem}\label{lemmeBiais}
Under conditions $(\ref{Ordre2})$ and  $(\ref{vng})$, we have $\sqrt{v_n} ( \bE( \gamtilde )  - \gF ) \ \tqdninf \ \lambda m$.
\end{lem} 
\zdeux \zdeux 

Note that Lemma \ref{lemmeDeltan} is a direct corollary of relation $(\ref{DecompDeltanm1})$ and of Lemmas  \ref{lemmeTLClambdan} and \ref{lemmeWnWtildenTLC}. Lemma \ref{lemmeTLClambdan} is included in the proof of Theorem 1 in \cite{GardesStupfler15}. We will provide the proofs of the other lemmas in the next subsections.
\zcinq

Let us now turn to the proof of Theorem \ref{TLCgamchap}. We have, thanks to Lemmas \ref{lemmegamchapgamtilde} and \ref{lemmeDeltan},
\begin{equation} \label{DecompGamChap1}
 \sqrt{v_n}(\gamchap - \gF) \; = \; \sqrt{v_n} (\Delta_n^{-1}\gamtilde - \gF) + \op1  \;
  =  \;  \Delta_n^{-1} \sqrt{v_n}\left( \, (\gamtilde - \gF) - \gF(\Delta_n-1)\,\right) + \op1 . \end{equation}
We consider
\[
 \Delta_n - 1  =  \frac{ \Fbarn(t_n) - \Fbartn }{ \Fbartn }    =  - \frac{ \Fn(t_n) -  F(t_n) }{ \Fbartn } 
\]
and we want to deal with this difference by introducing cumulative hazard functions (defined at the beginning of this section). But if there exists some data value $X_i^*$ which is both greater than $t_n$ and such that $nC_n(\Xiet) =1$, then $\Fn(t_n)=0$ and $\Lambda_n(t_n)$ is undefined. In order to avoid this, we introduce the variable 
$$ 
T= \max\{\Xiet ; nC_n(\Xiet) =1\}
$$ 
for which \cite{Stute2010} proved that $\bP(T=\min_{i\leq n} \Xiet)$ converges to $1$. Therefore, if we set $A_n=\{T\leq t_n\}$, then on $A_n$ we have $\Fn(t_n)>0$ on one hand, and on the other hand $\bP(A_n^c)\leq \bP(T\neq \min_{i\leq n} \Xiet) + \bP(\min_{i\leq n} \Xiet >t_n)$, which tends to $0$. We can thus write, using the mean value theorem, 
\begin{eqnarray*}
 \Delta_n - 1  & = &
 -\frac{ \exp(-\Lambda_n(t_n)) - \exp(-\Lambda(t_n)) } { \Fbartn }  \bI_{A_n} + \frac{ \Fbarn(t_n) - \Fbartn }{ \Fbartn } \bI_{A_n^c}
 \\ 
 & = & \xi_n \  \bI_{A_n} \frac{ \Lambda_n(t_n) - \Lambda(t_n) } { \Fbartn } \ + \ \frac{ \Fbarn(t_n) - \Fbartn }{ \Fbartn } \bI_{A_n^c} 
\end{eqnarray*}
where $\xi_n$ converges to 1 in probability, since both $\Lambda_n(t_n)$ and $\Lambda(t_n)$ converge to $0$. 
Therefore, using successively $\bP(A_n^c)\rightarrow 0$ and Lemmas \ref{lemmeContigLambdas}, \ref{lemmeTLClambdan} and \ref{lemmeWnWtildenTLC}, we can write
\begin{eqnarray}
 \sqrt{v_n} (\Delta_n - 1) & = & \xi_n \bI_{A_n} \sqrt{v_n}  \frac{ \hat\Lambda_n(t_n) - \Lambda(t_n) } { \Fbartn } + \op1 
 \; = \; \xi_n \bI_{A_n}  \sqrt{v_n} \overline{W}'_n + \op1 \nonumber \\
& = & \sqrt{v_n} \overline{W}'_n + \op1 .
\label{DecompDeltanm1}
\end{eqnarray}
On the other hand,  
\[
 \sqrt{v_n} (\gamtilde - \gF) = \sqrt{v_n} \overline{W}_n + \sqrt{v_n} (\bE(\gamtilde)-\gF) 
\]
and consequently, combining relations (\ref{DecompGamChap1}) and (\ref{DecompDeltanm1}) with Lemmas \ref{lemmeWnWtildenTLC} and \ref{lemmeBiais}, the theorem is proved :
\[
 \sqrt{v_n}(\gamchap-\gF) \; = \; \Delta_n^{-1} \left\{ \, \sqrt{v_n} (\overline{W}_n - \gF\overline{W}'_n) +  \sqrt{v_n} (\bE(\gamtilde)-\gF) + \op1 \right\} + \op1 \;  \stackrel{\cal L}{\longrightarrow} \; {\cal N} \left(\lambda m,s^2\right) .
\]

\subsection{Proof of Theorem \ref{TLCchapx}} 
Recall that $d_n=\frac{\Fb(t_n)}{p_n}\tinf$, and the notations $\Delta_n = \frac{\Fbarn(t_n)}{\Fbartn}$ (which satisfies (\ref{DecompDeltanm1}))  and   $v_n=n\Hb(t_n)$. We write
\[
 \frac{\chapx}{\xp} - 1   \; = \;  \frac{t_n}{\xp} (\Delta_n d_n)^{\gamchap} -1   
 \; = \; \Delta_n^{\gamchap} \left(   \frac{t_n}{\xp}\  d_n^{\gF} T_n^1  +  T_n^2 + T_n^3 \right), 
\zun
\]
where $T_n^1 := d_n^{\gamchap-\gF} -1$,  \ $T_n^2 :=  \frac{t_n}{\xp} d_n^{\gF} -1$ \ and \ $T_n^3 := 1 - \Delta_n^{-\gamchap}$. 
We are going to prove that both $T_n^2$ and $T_n^3$ are $o_{\bP}(\log d_n/\sqrt{v_n})$, and that $\frac{\sqrt{v_n}}{\log d_n} T_n^1  \stackrel{\cal L}{\longrightarrow}  {\cal N} \left(\lambda m , s^2 \right)$. This will conclude the proof, since both $\Delta_n$  and $ \frac{t_n}{\xp}\  d_n^{\gF}$ tend to $1$. 

Let us first focus on $T_n^1$. The mean value theorem yields 
\[
\frac{\sqrt{v_n}}{\log d_n} T_n^1 = \sqrt{v_n} (\gamchap - \gF) \exp(E_n), 
\]
where $| E_n |  \leq |\gamchap - \gF| \log d_n$ and therefore  $E_n$ tends to $0$ thanks to Theorem \ref{TLCgamchap} and  assumption $(\ref{condn})$. The desired result for $T_n^1$ is then implied by Theorem \ref{TLCgamchap}. 
\zun 

We now deal with $T_n^2$. Recalling that $\Fb(x)= x^{-\gF} l_F(x)$, by definition of $\xp$ we have
\[
T_n^2 = \left( \frac{l_F(\xp)}{l_F(t_n)} \right)^{-\gF} - 1
\]
We use  the following representation of $l_F$ (see \cite{Smith87} page 1195)  when $\rho <0$ : 
\[
l_F(x)= C \left(  1 + \rF^{-1} g(x) + o(g(x)) \right), \mbox{ for } x \rightarrow +\infty.
\]
Hence 
\[
\frac{l_F(\xp)}{l_F(t_n)} = 1  - \rF^{-1} g(t_n) \left( 1- \frac{g(\xp)}{g(t_n)}  + \op1 + o\left( \frac{g(\xp)}{g(t_n)}  \right)  \right). 
\]
But $g(\xp)/g(t_n)$ tends to $0$ because $\xp/t_n$ tends to infinity and 
\[
\left| \, g(\xp)/g(t_n) \, - \, \left( \xp/t_n \right)^{\rF} \, \right| \ \leq \ \sup_{y \geq 1}  \left| g(y t_n)/g(t_n) - y^{-\rho}\right| \ \longrightarrow \ 0. 
\]
It follows that $\frac{l_F(\xp)}{l_F(t_n)} = 1 - \rF^{-1} g(t_n) (1 + \op1) $. Thus  $\left| \left( l_F(\xp)/l_F(t_n) \right)^{-\gF} - 1 \right| \leq c  \left|l_F(\xp)/l_F(t_n)- 1 \right|$, for some constant $c$ and then 
\[
\frac{\sqrt{v_n}}{\log d_n}  |T_n^2| \leq c \rF^{-1} \sqrt{v_n}g(t_n)  \frac{1+\op1}{\log d_n}. 
\]
Assumption $(\ref{vng})$ and the fact that  $\log d_n$ tends to $0$ conclude the proof for $T_n^2$. 
\zdeux

For $T_n^3$, we use the mean value theorem to write 
\[
T_n^3 = \gamchap K_n^{-\gamchap-1} (\Delta_n -1),
\]
with $K_n$ tending to $1$. In view of (\ref{DecompDeltanm1}) and Lemma  \ref{lemmeWnWtildenTLC} , we thus have $\frac{\sqrt{v_n}}{\log d_n}  (\Delta_n -1) = O_{\bP} (1) / \log d_n = \op1$ and then the desired neglibility of $T_n^3$ follows.

\subsection{Proof of Lemma  \ref{lemmegamchapgamtilde}}

We have $\displaystyle
\Delta_n \gamchap =  \gamtilde + S_{n,1} + S_{n,2} , 
$
with 
\[
S_{n,1} :=  \unsurFbartn \moyennesurn \frac{\Fn (\Xiet) - F (\Xiet)}{C_n(\Xiet)} \  \logXiet \indicXiet 
\]
and 
\[
S_{n,2} :=  \unsurFbartn \moyennesurn F (\Xiet) \left( \frac{1}{C_n(\Xiet)}  - \frac{1}{C(\Xiet)} \right)  \  \logXiet \indicXiet.  
\]
Let us show that both $\sqrt{v_n} S_{n,1}$ and $\sqrt{v_n} S_{n,2}$ are $\op1$. 
On one hand, 
\begeq{majoSn1}
|\sqrt{v_n} S_{n,1}|  \leq    \left( \sqrt{n} \sup_{x > t_n} | \Fn(x) - F(x) | \right)  \ \sup_{\Xiet >t_n} \frac{C(\Xiet)}{C_n(\Xiet)}  \ \sqrt{\Hb (t_n)}  \bar{V}_n^1 
\fineq
where $\bar{V}_n^1 :=  \moyennesurn V_{i,n}^1$ with 
\[
V_{i,n}^1 := \unsurFbartn  \frac{\indicXiet}{C(\Xiet)} \ \logXiet. 
\]
Using (\ref{relationFetF}) and (\ref{defC}) yields 
\[
\bE(V_{i,n}^1) =  \unsurFbartn  \inttninf \frac{1}{F(x)} \log(x/t_n) dF(x) = (1 + \op1) \unsurFbartn  \inttninf \log(x/t_n) dF(x), 
\]
which converges to $\gF$ ; Markov inequality then yields $\sqrt{\Hb (t_n)}  \bar{V}_n^1 = \op1$. 
On the other hand, 
\begeq{majoSn2}
|\sqrt{v_n} S_{n,2}|  \leq  \sup_{ \Xiet >t_n} \frac{C(\Xiet)}{C_n(\Xiet)}  \ \left(  \sqrt{n} \sup_{ \Xiet >t_n}  | C_n(\Xiet) - C(\Xiet) | \right)  \   \sqrt{\Hb(t_n)}  \bar{V}_n^2, 
\fineq
where $\bar{V}_n^2 :=  \moyennesurn V_{i,n}^2$ with 
\[
V_{i,n}^2 := \unsurFbartn  \frac{F(\Xiet)}{C^2(\Xiet)} \ \logXiet  \indicXiet. 
\]
Using again (\ref{relationFetF}) and (\ref{defC}), we have
\[
\bE(V_{i,n}^2) =  p \unsurFbartn  \inttninf \frac{ \logxtntext}{F(x) \Gb(x)} dF(x) = p (1 + \op1) \unsurFbartn  \inttninf \frac{1}{\Gb (x)} \log(x/t_n) dF(x).  
\]
By Lemma \ref{lemmeIntlog} (where constant $c_1$ is defined), it comes $\bE(V_{i,n}^2) = (1+ \op1) \frac{pc_1}{\overline{G}(t_n)}$ and Markov inequality then yields {$\sqrt{\Hb(t_n)}  \bar{V}_n^2 =$} {$ O_{\bP}\left(({\Fbartn/\overline G(t_n)})^{1/2}\right)=\op1$}. 
Combining (\ref{majoSn1}) and (\ref{majoSn2}) with  Lemma \ref{lemmeFnC} ends the proof. 
\zdeux 
\zdeux 

\subsection{Proof of Lemma \ref{lemmeContigLambdas}}

Recall that  $T= \max\{\Xiet ; nC_n(\Xiet) =1\}$ and that we previously saw that $\bP(A_n)\rightarrow 1$ when $A_n=\{T\leq t_n\}$. Using the fact that $0\leq -\log(1-x)-x \leq \frac{x^2}{1-x}$  for any $0 \leq x < 1$, and that, on $A_n$, we have $nC_n(\Xiet) \geq 2$ for every $\Xiet > t_n$, we can write that 
\[
\frac{\sqrt{v_n}}{\Fbartn}  \bI_{A_n}  \ |\Lambda_n (t_n)- \hat{\Lambda}_n (t_n)| \; \leq \;
\frac{\sqrt{v_n}}{\Fbartn} \bI_{A_n}  \sum_{\Xiet > t_n} \frac{1}{n^2 C_n^2(\Xiet)} \  \frac{1}{1-\frac{1}{nC_n(\Xiet)}}  \; \leq \; 2  \bI_{A_n}  \frac{\sqrt{v_n}}{\Fbartn} \ \sum_{\Xiet > t_n} \frac{1}{n^2 C_n^2(\Xiet)} 
\]
Using Lemma \ref{lemmeFnC}, we have 
\[
\frac{\sqrt{v_n}}{\Fbartn} \ \sum_{\Xiet > t_n} \frac{1}{n^2 C_n^2(\Xiet)} \leq O_{\bP} (1) \sqrt{\frac{\Gb(t_n)}{n \Fbartn }} \moyennesurn \frac{\bI_{\Xiet > t_n}}{C^2(\Xiet)}. 
\]
Noting $Z_n = \moyennesurn \bI_{\Xiet > t_n} / C^2(\Xiet)$, and using (\ref{relationFetF}) and (\ref{defC}), we have
\[  
\bE(Z_n)=\inttninf  \frac{p}{F^2(x)} \frac{dF(x)}{\overline{G}(x)} = p (1 + \op1) \inttninf   \frac{dF(x)}{\Gb(x)} . 
\]
Via Lemma \ref{lemmeIntlog}, $\bE\left(\sqrt{\frac{\Gb(t_n)}{n \Fbartn }} Z_n\right)$ tends to $0$ and  therefore $\sqrt{\frac{\Gb(t_n)}{n \Fbartn }} Z_n = \op1$ by Markov's inequality, which ends the proof of the lemma. 
\zdeux \zdeux

\subsection{Proof of Lemma \ref{lemmeWnWtildenTLC}}

For brevity, we only prove the third part of the lemma. First, using relation (\ref{relationFetF}) and Lemma \ref{lemmeIntlog} (wherein the constants $c_0=1/q$, $c_1=\gF/q^2$, $c_2=2\gF^2/q^3$ are defined, with $q=1-\gF/\gG$), it is easily seen that 
\begin{eqnarray*}
\bE(V_{1,n}) & = & \unsurFbartn \inttninf \logxtn dF(x) 
\ \tqdninf \ \gF 
\\
\bE(V^2_{1,n}) & = & \frac{p}{\Fbartn^2} \inttninf \logxtndeux \frac{dF(x)}{\Gb(x)} 
\ = \ \frac{pc_2}{\Hb(t_n)} (1+o(1)) 
\\
\bE(V'_{1,n}) & = & \unsurFbartn \inttninf \frac{dF(x)}{F(x)} \ = \ \frac{\Lambda(t_n)}{\Fbartn}  
\ \tqdninf  \ 1 
\\ 
\bE((V'_{1,n})^2) & = & \unsurFbartndeux \inttninf \frac{p}{\Gb(x)F^2(x)} dF(x) \ = \ 
\frac{p(1+o(1))}{\Fbartn^2} \inttninf \frac{dF(x)}{\Gb(x)} 
\ = \ \frac{pc_0}{\Hb(t_n)} (1+o(1)) 
\\
\bE(V_{1,n}V'_{1,n}) & = & 
\frac{p(1+o(1))}{\Fbartn^2} \inttninf \logxtn \frac{dF(x)}{\Gb(x)} 
\ = \ \frac{pc_1}{\Hb(t_n)} (1+o(1)) 
\end{eqnarray*}
Introducing $U_{i,n}=W_{i,n}-\gF W'_{i,n}$ and  $S_n=\sum_{i\leq n} U_{i,n}$, we thus obtain ($s^2$ is defined in the statement of the lemma)
\[
 {\bV}ar (U_{1,n}) = {\bV}ar(V_{1,n}-\gF V'_{1,n}) = \frac{s^2}{\Hb(t_n)}(1+o(1)) 
 \]
and consequently $\sqrt{v_n}(\overline W_n - \gF \overline W'_n)= \sqrt{v_n}  S_n/n = s(1+o(1)) S_n / {\bV}ar(S_n)$, which converges in distribution to ${\cal N}(0,s^2)$ as soon as Lyapunov's condition holds.  
After some simplifications,  Lyapunov's condition becomes the existence of some $\delta>0$ such that \[
 n^{-\delta/2} (\Hb(t_n))^{1+\delta/2}\bE(|U_{1,n}|^{2+\delta})\;\tqdninf \; 0. 
\]
Proceeding as in \cite{GardesStupfler15}, and noting that $\bE(V_{1,n})-\gF\bE(V'_{1,n})$ vanishes to $0$, the double application of the inequality $|a+b|^{2+\delta}\leq 2^{1+\delta}(|a|^{2+\delta}+|b|^{2+\delta})$  shows that it suffices to prove the following, for some $\delta>0$ : 
\begeq{condLyap}
  n^{-\delta/2} (\Hb(t_n))^{1+\delta/2}\bE(|V|^{2+\delta}) \ \tqdninf \ 0
  \hspace{0.3cm} \mbox{for both $V=V_{1,n}$ and $V=V'_{1,n}$}
\fineq
We prove this property for $V=V_{1,n}$, the proof for $V=V'_{1,n}$ being very similar. We have
\[
 \bE(|V_{1,n}|^{2+\delta}) \; = \; p^{2+\delta} (\Fbartn)^{-2-\delta} \inttninf \log^{2+\delta}\left(\frac x{t_n}\right) \frac{dF(x)}{\Gb^{1+\delta}(x)} 
\]
Mimicking the proof of Lemma \ref{lemmeIntlog} stated in the appendix, and because $\delta$ can be chosen arbitrary small (so that $(1+\delta)/\gG$ remains lower than $1/\gF$), we can prove that 
\[
 \frac{\Gb^{1+\delta}(t_n)}{\Fbartn} \inttninf \log^{2+\delta}\left(\frac x{t_n}\right) \frac{dF(x)}{\Gb^{1+\delta}(x)}  \; = \; O(1)
\]
and therefore, since we assumed that $n\Hb(t_n)\tinf$, the desired property (\ref{condLyap}) holds for  $V=V_{1,n}$ : 
\[
 n^{-\delta/2} (\Hb(t_n))^{1+\delta/2}\bE(|V_{1,n}|^{2+\delta}) 
 \ \leq \ O(1)  n^{-\delta/2} (\Hb(t_n))^{1+\delta/2} (\Fbartn)^{-2-\delta}  \Fbartn \Gb^{-1-\delta}(t_n)
 \ = \ O(1) (n\Hb(t_n))^{-\delta/2} \ \tqdninf \ 0.
\]

\subsection{Proof of Lemma \ref{lemmeBiais}}

Recall that $\bE(\gamtilde)= \unsurFbartn \inttninf \logxtn dF(x) = \int_1^{+\infty} \frac 1 y \ \frac{\Fb(y t_n)}{\Fbartn}dy$ by integration by parts and change of variables. Since $\Fb(y) = y^{-1/\gF} l_F(y)$, we have 
\[ 
\sqrt{v_n} (\bE(\gamtilde)- \gF) =  \sqrt{v_n} \int_1^{+\infty} y^{-1/\gF -1} \left(  \frac{l_F(y t_n)}{l_F(t_n)} -1 \right) \ dy, 
             \] 
and using  assumption $(\ref{Ordre2})$  and  Proposition 3.1 in \cite{Smith87}, we can write 
\[
\int_1^{+\infty} y^{-1/\gF -1} \left(  \frac{l_F(y t_n)}{l_F(t_n)} -1 \right) \ dy = g(t_n) \  \int_1^{+\infty} y^{-1/\gF -1} h_{\rF} (y) dy + o(g(t_n)). 
\]
The result then follows from assumption $(\ref{vng})$  and the fact that $ \int_1^{+\infty} y^{-1/\gF -1} h_{\rF} (y) dy \,=\, m$.

\section{Appendix}

This appendix contains two lemmas : Lemma \ref{lemmeFnC} contains results which are proved elsewhere but are crucial for our proof, and which we thus restate here, whereas Lemma \ref{lemmeIntlog} is  a variant of a particular case of Lemma 2 in \cite{GardesStupfler15}, and states essential equivalences for our proofs.
 \zdeux

\begin{lem}\label{lemmeFnC}  
If $t_n$ tends to infinity with $n$, then
\begitem
\item[(a)] $\sqrt{n} \sup_{x > t_n} | \Fn(x) - F(x) | = O_{\bP}(1)$.
\item[(b)] $\sup_{1 \leq i \leq n} \bigg\{ \frac{C(\Xiet)}{C_n(\Xiet)} \, \bigg| \, \Xiet >t_n \bigg\}  = O_{\bP}(1)$.
\item[(c)] $\sqrt{n} \sup_{1 \leq i \leq n}  \left\{ \; | C_n(\Xiet) - C(\Xiet) | \; \big| \; \Xiet >t_n \; \right\} = O_{\bP}(1)$. 
\finit 
\end{lem}
{\it Proof}
\zun\\
$(a)$ is a consequence of  point 6 page 176  in \cite{Woodroofe85}. $(b)$ is proved in \cite{GardesStupfler15} (see lemma 5), following the ideas contained in \cite{Stute93}. Since $C_n= F_n^*  - G_n^*$ , where  $F_n^*$ and $G_n^*$ are respectively the empirical distribution functions of $F^*$ and $G^*$, $(c)$ is a consequence of $\sqrt{n} \sup_{x \geq 0}  | F_n^*(x) - F^*(x) |=O_{\bP}(1)$ and $\sqrt{n} \sup_{x \geq 0}  | G_n^*(x) - G^*(x) |=O_{\bP}(1)$ (see \cite{Woodroofe85} pages 172-173). 

\begin{lem}\label{lemmeIntlog} 
Under condition $(\ref{Ordre1})$, for  any $k\in\bN$, as $n\tinf$,
$$
 \inttninf \log^k \left(\frac x {t_n}\right) \, \frac{dF(x)}{\Gb (x)}   =  c_k \frac{\Fbartn}{\Gb (t_n)} (1 + o(1))
$$
where $\displaystyle c_k=\frac{\gF^k\,k!}{(1-\gF/\gG)^{k+1}}$. 
\end{lem}
{\it Proof}
\zun\\
Let us note $\alpha=1/\gG$ and $\beta=1/\gF$, which satisfy $0<\alpha<\beta$ by assumption. We need to prove that the following quantity converges to $c_k$ (below, $\delta >0$ is arbitrary small)
\begin{eqnarray}
& & \hspace*{-2.cm}\frac{\overline G(t_n)}{\overline F(t_n)} \inttninf \log^k\left(\frac x{t_n}\right) \frac{dF(x)}{\overline G(x)}
\nonumber\\
& = &
-\intuninf \log^k(y) \frac{\overline G(t_n)}{\overline G(yt_n)} \frac{t_n d\overline F(yt_n)}{\Fbartn}
\nonumber \\
& = &
-\intuninf \log^k(y) \; y^{\alpha} \; \frac{t_n d\overline F(yt_n)}{\Fbartn}
\nonumber \\
& & \hspace*{1.cm} -\intuninf \log^k(y) y^{\alpha+\delta} \left\{ \frac{\overline G(t_n)}{\overline G(yt_n)} \frac{(yt_n)^{-\alpha-\delta}}{t_n^{-\alpha-\delta}} - y^{-\delta} \right\}
\frac{t_n d\overline F(yt_n)}{\Fbartn}
\nonumber \\
& = & I_{n,k}(\alpha) + o(1)I_{n,k}(\alpha+\delta)
\label{relationInk}
\end{eqnarray}
In the last line, we used Theorem 1.5.2 in \cite{BinghamGoldieTeugels} with the fact that $x\mapsto x^{-\alpha-\delta}/\overline G(x)$ is regularly varying of order $-\delta$. It thus remains to prove that $I_{n,k}(\alpha)$ converges to $c_k$ (the same being true for $I_{n,k}(\alpha+\delta)$).
We now introduce the notations : for $\theta>0$
\[
J_k(\theta) = \intuninf \log^k(y)y^{-\theta-1}dy = \frac{k!}{\theta^{k+1}} \makebox[2.cm][c]{and}
J_{n,k} = \intuninf \log^k(y)y^{\alpha-1}\frac{\overline F(yt_n)}{\Fbartn} dy .
\]
For any $\delta\in]0,\beta-\alpha[$, since the function $x\mapsto x^{\beta-\delta}\overline F(x)$ is regularly varying of order $-\delta$, we have
\begin{eqnarray*}
J_{n,k} & = & \intuninf \log^k(y) y^{\alpha-\beta-1}dy \ + \
\intuninf \log^k(y) y^{\alpha-1} \left( \frac{\overline F(yt_n)}{\Fbartn}\frac{(yt_n)^{\beta-\delta}}{t_n^{\beta-\delta}} - y^{-\delta} \right) y^{-\beta+\delta}dy \\
& = & J_k(\beta-\alpha) + o(1)
\end{eqnarray*}
We thus have, by integration by parts and the relation $kJ_{k-1}(\theta) = \theta J_k(\theta)$,
\begin{eqnarray*}
I_{n,k}(\alpha) & = & \intuninf (k\log^{k-1}(y) + \alpha \log^k(y)) y^{\alpha-1} \frac{\overline F(yt_n)}{\Fbartn} dy
\\
& = & k\, J_{n,k-1} + \alpha\, J_{n,k}
\\
& \tqdninf & kJ_{k-1}(\beta-\alpha) + \alpha J_k(\beta-\alpha) \ = \ \beta J_k(\beta-\alpha)
\ = \ \frac 1 {\gF} \frac{k!}{(\gF^{-1}-\gG^{-1})^{k+1}} \ = \ c_k
\end{eqnarray*}

\bibliographystyle{elsarticle-num}
\bibliography{Bibtruncation.bib}{}

\begin{thebibliography}{10}
\expandafter\ifx\csname url\endcsname\relax
  \def\url#1{\texttt{#1}}\fi
\expandafter\ifx\csname urlprefix\endcsname\relax\def\urlprefix{URL }\fi
\expandafter\ifx\csname href\endcsname\relax
  \def\href#1#2{#2} \def\path#1{#1}\fi

\bibitem{BeirlantDierckxFVGuillou07}
J.~Beirlant, G.~Dierckx, A.~Fils-Villetard, A.~Guillou, Estimation of the
  extreme value index and extreme quantiles under random censoring, Extremes 10
  (2007) 151--174.

\bibitem{NecirTroncatureHill15}
S.~Benchaira, D.~Meraghni, A.Necir, On the estimation of the extreme value
  index for randomly right-truncated data and application.

\bibitem{BinghamGoldieTeugels}
N.~Bingham, C.~Goldie, J.~Teugels, Regular variation, Cambridge University
  Press, 1987.

\bibitem{EinmahlFVGuillou08}
J.~Einmahl, A.~Fils-Villetard, A.~Guillou, Statistics of extremes under random
  censoring, Bernoulli 14 (2008) 207--227.

\bibitem{GardesStupfler15}
L.~Gardes, G.~Stupfler, Estimating extreme quantiles under random truncation,
  TEST 24 (2015) 207--227.

\bibitem{deHaanFerreira2006}
L.~de~Haan, A.~Ferreira, Extreme Value Theory : an introduction, Springer
  Series in Operations Research and Financial Engineering, Springer, 2006.

\bibitem{Smith87}
R.~Smith, Estimating tails of probability distributions, Annals of Statistics
  15~(3) (1987) 1174--120.

\bibitem{Stute2010}
E.~Strzalkowska-Kominiak, W.~Stute, On the probability of holes in truncated
  samples, Journal of Statistical Planning and Inference 140 (2010) 1519--1528.

\bibitem{Stute93}
W.Stute, Almost sure representations of the product-limit estimator for
  truncated data, Annals of statistics 21~(1) (1993) 146--156.

\bibitem{Woodroofe85}
M.~Woodroofe, Estimating a distribution function with truncated data, Annals of
  statistics 13~(1) (1985) 163--177.

\bibitem{Worms14}
J.~Worms, R.~Worms, New estimators of the extreme value index under random
  right censoring, for heavy-tailed distributions, Extremes 17 (2014) 337--358.

\end{thebibliography}

\end{document}